\documentclass[a4paper,pdftex,reqno,10pt]{amsart}
\usepackage{geometry}
\usepackage{graphicx}
\usepackage{multirow}
\usepackage{times}
\usepackage{amssymb,amsmath}

\def\Ps{\mathbb{P}}

\def\@begintheorem#1#2{\list{}{\thm@body}%
  \item[]{\bf #1~#2.}\quad\it\ignorespaces}
\def\@opargbegintheorem#1#2#3{\list{}{\thm@body}%
  \item[]{\bf #1~#2~\ifrembrks #3\global\rembrksfalse\else (#3)\fi.}%
  \quad\it\ignorespaces}
\def\@endtheorem{\endlist}
\newtheorem{theorem}{Theorem}[section]

\newtheorem{corollary}[theorem]{Corollary}
\newtheorem{lemma}[theorem]{Lemma}

\newtheorem{definition}[theorem]{Definition}

\newtheorem{conjecture}[theorem]{Conjecture}
\newcommand{\eop}{\hfill{$\Box$}}

\newcommand{\OZ}{OZ}

\providecommand{\Aut}{\mathop{\rm Aut}\nolimits}

\begin{document}

\title{Integral point sets in higher dimensional affine spaces over finite fields}
\author{Sascha Kurz}
\address{University of Bayreuth, Department of Mathematics, D-95440 Bayreuth, Germany, sascha.kurz@uni-bayreuth.de, www.wm.uni-bayreuth.de/index.php?id=sascha}
\author{Harald Meyer}
\address{University of Bayreuth, Department of Mathematics, D-95440 Bayreuth, Germany}
\begin{abstract}
    We consider point sets in the $m$-dimensional affine space $\mathbb{F}_q^m$
    where each squared Euclidean distance of two points is a square in $\mathbb{F}_q$. It turns out that the situation in
    $\mathbb{F}_q^m$ is rather similar to the one of integral distances in Euclidean spaces. Therefore we expect the
    results over finite fields to be useful for the Euclidean case.

    We completely determine the automorphism group of these spaces which preserves
    integral distances. For some small parameters $m$ and $q$ we determine the maximum cardinality
    $\mathcal{I}(m,q)$ of integral point sets in $\mathbb{F}_q^m$. We provide upper bounds
    and lower bounds on $\mathcal{I}(m,q)$. If we map integral distances
    to edges in a graph, we can define a graph $\mathfrak{G}_{m,q}$ with vertex set $\mathbb{F}_q^m$. It turns out that
    $\mathfrak{G}_{m,q}$ is strongly regular for some cases.
    
    \medskip
    
    \noindent
    \textbf{Keywors:} finite geometry, integral distances, integral point sets, automorphism group, strongly regular graphs
    
\end{abstract}

\maketitle

\section{Introduction and Notation}
\noindent
Integral point sets, i.~e.{} point sets with pairwise integral distances in Euclidean space, have been considered since
the time of the Pythagoreans, who studied rectangles with integral side lengths and integral diagonal. Even nowadays there are a lot of unsolved problems concerning integral point sets \cite[section 5.11]{1086.52001}. E.~g.{} it is not known whether a perfect cuboid, that is a box with integer edges, face diagonals, and body diagonal, exists \cite[Problem D18]{UPIN}.

Applications originate from chemistry (molecules), physics (wave lengths), robotics, architecture, see
\cite{integral_distances_in_point_sets}. The concept of integral point sets can be generalized to commutative rings in order to study the underlying structure suppressing number theoretical difficulties. Here, we consider point sets in $m$-dimensional affine space $\mathbb{F}_q^m$ where each squared Euclidean distance of two points is a square in $\mathbb{F}_q$. It turns out that the situation in $\mathbb{F}_q^m$ is rather similar to the one in Euclidean spaces. Therefore we expect the results over finite fields to be useful for the Euclidean case.

In this context we would like to remark the famous problem of P. Erd\H{o}s who asked for seven points in the plane, no three on a line, no four on a circle with pairwise integral distances \cite[Problem D20]{UPIN}, \cite{1062.00002}. Several conjectures and incorrect proofs circulated that such a point set cannot exist. In \cite{axel_2} the authors find corresponding examples over $\mathbb{F}_q^2$ consisting of nine points which finally lead to the discovery of an integral heptagon in the Euclidean plane
in \cite{kreisel}.

There has been done a lot of work on integral point sets in Euclidean spaces, see e.~g.{} \cite{integral_distances_in_point_sets,kreisel,phd_kurz,paper_laue,sascha-alfred}. Some authors also consider other spaces, e.~g.{} Banach spaces \cite{banach}, integral point sets over rings \cite{axel_2}, or integral point sets over finite fields \cite{algo,Dimiev-Setting,michael_1,integral_over_fields}. In \cite{integral_over_fields} one of the authors of this article determines the automorphism group for dimension $m=2$, and in \cite{michael_1} integral point sets over $\mathbb{F}_q^2$, which are maximal with respect to inclusion, were classified for $q\le 47$. For $m=2$ and $q\equiv 3\pmod 4$ the graphs $\mathfrak{G}_{m,q}$ from Section \ref{sec_graph_of_integral_distances} are isomorphic to Paley graphs of square order. So in some sense these graphs $\mathfrak{G}_{m,q}$ generalize Paley graphs. In \cite{Blokhuis-1984} Blokhuis has determined the structure of cliques of maximal size in Paley graphs of square
order. Since, whenever two points in $\mathbb{F}_q^m$ are at integral distance the whole line through these points is an integral point set, integral point sets over $\mathbb{F}_q^m$ correspond to point sets with few directions being contained in a given set. From this point of view methods from R\'edei's seminal work \cite{0229.12019} can be applied to obtain results for point sets with few directions, see e.~g.{} \cite{0561.12009,0945.51002,1117.05018}.

Here, after giving the basic facts on integral point sets over affine planes in Section \ref{sec_integral_point_sets}, we completely determine the automorphism group of $\mathbb{F}_q^m$ with respect to integral distances in Theorem \ref{thm_automorphismgroup_general} and analyze its operation on $\mathbb{F}_q^m$ in Section \ref{sec_automorphism_group}. We introduce and analyze the graphs of integral distances $\mathfrak{G}_{m,q}$ for $m\ge 3$, $2\nmid q$ in Section \ref{sec_graph_of_integral_distances}. They arise from $3$-class association schemes. The determination of some of their parameters let us conjecture that they are strongly regular for even dimensions $m$. In Section \ref{sec_maximum_cardinality} we consider the maximum cardinality $\mathcal{I}(m,q)$ of integral point sets over $\mathbb{F}_q^m$ and provide some new exact numbers for dimension $m = 3$. For general dimension $m$ we state upper bounds and some constructions yielding lower bounds.
We finish with a conclusion and an outlook in Section \ref{sec_conclusion}.

We end this introduction with some notation we will keep throughout the paper. Let $p$ be a prime and let $q=p^r$ be a power. We write $\mathbb{F}_q$ for the field with $q$ elements and $\mathbb{F}_q^* := \mathbb{F}_q \backslash \{0\}$ for the units of $\mathbb{F}_q$. Our notation for the general linear group, i.~e.{} the set of invertible $m$ by $m$ matrices over $\mathbb{F}_q$, is $GL(m,q)$. By
\[
  O(m,q):=\left\{A\in GL(m,q) \mid A^TA=AA^T=E^{(m)}\right\},
\]
where $E^{(m)}$ is the $m\times m$ identity matrix, we denote the orthogonal group in dimension $m$. We remark that for even dimension $2n$ the orthogonal group comes in two types $O^+(2n,q)$ and $O^-(2n,q)$, and the group we defined above is isomorphic to $O^+(2n,q)$ in this case. By $A\Gamma L(m,q)$ we denote the affine general semilinear group over $\mathbb{F}_q$ and by
\[
  \OZ(m,q):=\left\{A\in GL(m,q) \mid A^TA=AA^T\in\mathbb{F}_q^*\cdot E^{(m)}\right\},
\]
we denote the smallest group containing $O(m,q)$ and the center of $GL(m,q)$, i.~e.{} the diagonal matrices with equal entries at the diagonal.\par

Originally integral point sets were defined in $m$-dimensional Euclidean spaces $\mathbb{E}^m$ as sets of $n$ points with pairwise integral distances in the Euclidean metric, see e.~g.{} \cite{integral_distances_in_point_sets,kreisel,phd_kurz,paper_laue,sascha-alfred} for an overview on the most recent results. Here we consider integral point sets in the affine spaces $\mathbb{F}_q^m$. We equip those spaces with a bilinear form
\[
  \langle u, v\rangle := u^Tv=\sum_{i=1}^m u_i v_i
\]
and a squared distance
\[
  d^2(u,v):= \langle u - v, u - v\rangle =(u-v)^T(u-v)= \sum_{i=1}^m (u_i-v_i)^2\quad\in\mathbb{F}_q
\]
for any two points $u=\begin{pmatrix}u_1&\dots&u_m\end{pmatrix}^T$, $v=\begin{pmatrix}v_1&\dots&v_m\end{pmatrix}^T$ in $\mathbb{F}_q^m$. We say that two points $u, v\in \mathbb{F}_q^m$ are at integral distance if $d^2(u,v)$ is contained in the set $\square_q:=\left\{\alpha^2\mid \alpha\in\mathbb{F}_q\right\}$ consisting of the squares in $\mathbb{F}_q$. As in the Euclidean space we define the cross product of two vectors $u, v\in\mathbb{F}_q^3$ by
\[
  u\times v := \begin{pmatrix}(u_2 v_3 - u_3 v_2)&(-u_1 v_3 + u_3 v_1)&(u_1 v_2 - u_2 v_1)\end{pmatrix}^T
  \in \mathbb{F}_q^3.
\]
The proofs of some of the common formulas for the cross product do not depend on any specific attributes of the Euclidean space, so they still hold in $\mathbb{F}_q^3$. Especially this is true for the formulas
\[
  \langle u\times v, u\rangle = \langle u\times v, v\rangle = 0
\]
and
\[
  \langle u\times v, u\times v\rangle = \langle u, u\rangle \cdot \langle v, v\rangle - \langle u, v\rangle^2
\]
we will use later on. The notation
\[
  U^\bot := \{ v\in \mathbb{F}_q^m \mid \langle u, v\rangle = 0 \mbox{ for all } u\in U\}
\]
for a subspace $U\subseteq\mathbb{F}_q^m$ is also inspired by the notation for the Euclidean space. As a shorthand
we use $u^\bot$ instead of $\{u\}^\bot$ for a vector $u\in\mathbb{F}_q^m$.\par

We have the equation
\[
  \langle Au, Av\rangle =(Au)^T(Av)=u^TA^TAv=u^Tv = \langle u, v\rangle
\]
for all $u,v\in\mathbb{F}_q^m$ and all $A\in O(m,q)$. If we have a matrix $A\in GL(m,q)$ with
$\langle Au, Av\rangle = \langle u, v\rangle$ for all $u, v\in\mathbb{F}_q^m$, then we have
\[
  u^T A^T A v = (Au)^T (Av) = \langle Au, Av \rangle = \langle u, v\rangle = u^T v
\]
for all $u, v\in\mathbb{F}_q^m$, i.~e.{} $A^T A = E^{(m)}$ and so $A$ is an element of $O(m,q)$.


\section{Integral point sets}
\label{sec_integral_point_sets}

\noindent
%
%

A set $\Ps$ of points in $\mathbb{F}_q^m$ is called an integral point set if all pairs of points are at integral distance, i.~e.{} if $d^2(u,v)\in\square_q$ for all $u,v\in\Ps$. As a shorthand we define
$\Delta:\mathbb{F}_q^m\times\mathbb{F}_q^m\rightarrow\{0,1\},$
\[
  (u,v)\mapsto \left\{\begin{array}{ll} 1 & \text{if $u$ and $v$ are at integral distance},\\
  0 & \text{otherwise}.\end{array}\right..
\]
By $\mathcal{I}(m,q)$ we denote the maximum cardinality of an integral point set in $\mathbb{F}_q^m$.




\begin{lemma}
  \label{lemma_first_bounds}
  \[
    q\le\mathcal{I}(m,q)\le q^m.
  \]
\end{lemma}
\textit{Proof.}
  For the lower bound we consider the \textit{line} $\Ps=\left\{\begin{pmatrix}\alpha&0&\dots&0\end{pmatrix}^T
  \mid \alpha\in\mathbb{F}_q\right\}$.
\hfill{$\square$}

\begin{lemma}
  If $2\mid q$ then we have
  $\mathcal{I}(m,q)=q^m$.
\end{lemma}
\textit{Proof.}
  For two points $u=\begin{pmatrix}u_1&\dots&u_m\end{pmatrix}^T$, $v=\begin{pmatrix}v_1&\dots&v_m\end{pmatrix}^T$
  in $\mathbb{F}_q^m$ we have
  \[
    d^2(u,v)=\sum_{i=1}^m \left(u_i-v_i\right)^2=\underset{\in\mathbb{F}_q}{\underbrace{\left(\sum_{i=1}^m
    u_i+v_i\right)}}^2.
  \]
\hfill{$\square$}

So in the remaining part of this article we consider only the cases where $2\nmid q$.

The lower bound of Lemma \ref{lemma_first_bounds} is attained in some cases, too, see e.~g.{} \cite{integral_over_fields} for a proof:
\begin{theorem}
  \label{thm_cardinality_m_2}
  \[
    \mathcal{I}(2,q)=q\text{ for }2\nmid q.
  \]
\end{theorem}

\section{Automorphisms preserving integral distances}
\label{sec_automorphism_group}

The primary object of this section is to determine the automorphism
group of $\mathbb{F}_q^m$ with respect to $\Delta$. We first have to 
define what we consider as an automorphism.
\begin{definition}
  \label{def_equivalence}
  An automorphism of $\mathbb{F}_q^m$ with respect to $\Delta$ is a 
bijective mapping
  $\varphi\in\text{A$\Gamma$L($\mathbb{F}_q$,m)}$
  with
  \[
    \Delta\left(u,v\right)=\Delta\left(\varphi(u),\varphi(v)\right)
  \]
  for all $u,v\in\mathbb{F}_q^m$. The group of automorphisms with 
respect to $\Delta$ is denoted by $\Aut(m,q)$.
\end{definition}

In other words this definition says that $\varphi$ has to map affine
subspaces, like points, lines, or hyperplanes, to affine subspaces with
equal dimension, and has to preserve the integral distance property.
It is easy to find the automorphisms with respect to $\Delta$, see Lemma
\ref{lemma_equivalence}. The difficult part is to show that these are all
automorphisms for $m\ge 3$. Here, the translations and Frobenius
homomorphisms do not cause any problems, so the most important part of
the determination of the automorphism group is to determine the linear
automorphisms. The main theorem of this section is:

\begin{theorem}
  \label{thm_automorphismgroup_general}
  \[
   \Aut(m,q)\cap GL(m,q) = \OZ(m,q) \text{ for }2\nmid q
   \text{ and }m\ge 3.
  \]
\end{theorem}

The idea of the proof of this theorem is to do induction
using the results for the case of dimension $m = 2$ one of the authors achieved in
\cite{integral_over_fields}. But for the proof we need a lot of facts about squares in
and about the action of the automorphism group on $\mathbb{F}_q^m$. We
will prove these facts step by step in several lemmas. We start with
some statements about roots in $\mathbb{F}_q$ and the set of solutions of
quadratic equations in $\mathbb{F}_q$.

\begin{definition}
  Let $q\equiv 1\mod 4$. By $\omega_q$ we denote an element with $\omega_q^2=-1$.
\end{definition}

\begin{lemma}
  \label{lemma_root_of_minus_one}
  For a finite field $\mathbb{F}_q$ with $q=p^r$ and $p\neq 2$ we have $-1\in\square_q$ iff $q\equiv 1\mod 4$,
  $\omega_q\in\square_q$ iff $q\equiv 1\mod 8$, and $2\in\square_q$ iff $q\equiv \pm 1\mod 8$.
\end{lemma}
\textit{Proof.}
  The multiplicative group of the units $\mathbb{F}_q^*$ is cyclic of order $q-1$. Elements of order $4$ are exactly
  those elements $\alpha$ with $\alpha^2=-1$. A similar argument holds for the fourth roots of $-1$. For the
  last statement
  we have to generalize the second auxiliary theorem (Erg\"anzungssatz) of the quadratic reciprocity law to
  $\mathbb{F}_q$: If $q = p$ is prime then the statement is true by the second auxiliary theorem. If $2$ is a square
  in $\mathbb{F}_p$ then $2$ is also a square in $\mathbb{F}_{p^r}$ and from $p\equiv \pm 1 \mod 8$ we get
  $q = p^r \equiv \pm 1 \mod 8$, i.~e.{} the statement is true in this case. If $2$ is not a
  square in $\mathbb{F}_p$ then the polynomial $x^2 - 2\in \mathbb{F}_p[x]$ is irreducible and $2$ is a square in
  $\mathbb{F}_p[x]/\left(x^2 - 2\right) \cong \mathbb{F}_{p^2}$. Hence $2$ is a square in $\mathbb{F}_{p^k}$ iff
  $2\mid k$. As $p$ is an odd prime with $p\equiv \pm 3 \mod 8$ we obtain $p^2 \equiv 1 \mod 8$ and also
  $p^{2k} \equiv 1 \mod 8$ as well as $p^{2k+1} \equiv \pm 3$ in this case.
\hfill{$\square$}

\begin{definition}
A triple $(\alpha,\beta,\gamma)$ is called \emph{Pythagorean triple over $\mathbb{F}_q$} if $\alpha^2+\beta^2=\gamma^2$.
\end{definition}

In the following it will be useful to have a parametric representation of the Pythagorean triples over $\mathbb{F}_q$.

\begin{lemma}
  \label{lemma_pythagorean_triples}
  For $2\nmid q$ let $\gamma\in\mathbb{F}_q$ and let $\mathbb{H}_\gamma$ be the set of Pythagorean triples
  $(\alpha,\beta,\gamma)$ over $\mathbb{F}_q$.
  \begin{enumerate}
  \item[(a)] If $\gamma=0$ then\\
  $\begin{array}{lll}
  \mathbb{H}_0 & = \begin{cases}\{(\tau,\pm \tau\omega_q,0)\mid \tau\in\mathbb{F}_q\} & \mbox{if }
  q\equiv 1\pmod 4,\\\{(0,0,0)\} &
  \mbox{if }q\equiv 3\pmod 4,\end{cases} \\
  \mbox{and} &&\\
  \lvert \mathbb{H}_0\rvert & = \begin{cases}2q-1 & \mbox{if }q\equiv 1\pmod 4, \\ 1 & \mbox{if }
  q\equiv 3\pmod 4.\end{cases}
  \end{array}$
  \item[(b)] If $\gamma\neq 0$ then\\
  $\begin{array}{lll}
  \mathbb{H}_\gamma & = \{(\pm \gamma,0,\gamma)\}\cup\{(0,\pm \gamma,\gamma)\}\cup
  \left\{\left(\frac{\tau^2-1}{\tau^2+1}\cdot \gamma,\frac{2\tau}{\tau^2+1}
  \cdot \gamma,\gamma\right)\mid\tau\in\mathbb{F}_q^*,\tau^2\neq \pm 1\right\}\\
  \mbox{and} &&\\
  \lvert \mathbb{H}_\gamma\rvert & = \begin{cases}q-1 & \mbox{if }q\equiv 1\pmod 4,\\q+1 & \mbox{if }
  q\equiv 3\pmod 4.\end{cases}\end{array}$
  \item[(c)]
  There are exactly $q^2$ Pythagorean triples over $\mathbb{F}_q$.
  \end{enumerate}
\end{lemma}
\textit{Proof.}
  Part~(a) is easy to verify.
  For part~(b) there are 4 solutions with $\alpha\beta=0$, these are $\{(0,\pm \gamma,\gamma),
  (\pm \gamma,0,\gamma)\}$.
  For $\alpha\beta\neq 0$ we get:
  \[
    \alpha^2+\beta^2=\gamma^2\quad\Leftrightarrow\quad\frac{\gamma-\alpha}{\beta}\cdot\frac{\gamma+\alpha}{\beta}=1.
  \]
  Setting $\tau:=\frac{\gamma+\alpha}{\beta}\in\mathbb{F}_q^*$ we obtain $\tau^{-1}=\frac{\gamma-\alpha}{\beta}$, hence
  \[
    \frac{\alpha}{\beta}=\frac{\tau-\tau^{-1}}{2}\quad\mbox{and}\quad\frac{\gamma}{\beta}=\frac{\tau+\tau^{-1}}{2}.
  \]
  Because of $\alpha\neq 0$, $\gamma\neq 0$ we have $ \tau \not= \pm \tau^{-1}$, i.~e.{} $\tau^2\notin \{-1,1\}$.
  It follows
  \[
    \alpha=\frac{\tau-\tau^{-1}}{\tau+\tau^{-1}}\cdot \gamma\quad\text{and}\quad
    \beta=\frac{2}{\tau+\tau^{-1}}\cdot \gamma.
  \]
  It is easily checked that for all admissible values of $\tau$, the resulting triples $(\alpha,\beta,\gamma)$
  are pairwise different Pythagorean triples.

  The expression for the number of solutions follows because $-1$ is a square in $\mathbb{F}_q$ exactly if
  $q\equiv 1\pmod 4$.

  With part~(a) and part~(b) we get the number of Pythagorean triples over $\mathbb{F}_q$ as
  \[
  \sum_{\gamma\in\mathbb{F}_q}\lvert \mathbb{H}_\gamma\rvert = \lvert \mathbb{H}_0\rvert +
  (q-1)\cdot\lvert \mathbb{H}_1\rvert = q^2.
  \]
  So also part~(c) is shown.
\hfill{$\square$}

From this lemma we can deduce the following corollary.
\begin{corollary}
  \label{cor_num_sol}
  If $\mathbb{I}_\gamma:=\left\{(\alpha,\beta)\mid \alpha^2+\beta^2=\gamma\right\}$ then we have
  \begin{align*}
  \lvert \mathbb{I}_0\rvert & = \begin{cases}2q-1 & \mbox{if }q\equiv 1\pmod 4, \\ 1 &
  \mbox{if }q\equiv 3\pmod 4,\end{cases}
  \end{align*}
  and
  \begin{align*}
  \lvert \mathbb{I}_\gamma\rvert & = \begin{cases}q-1 & \mbox{if }q\equiv 1\pmod 4,\\q+1
   & \mbox{if } q\equiv 3\pmod 4\end{cases}
  \end{align*}
  for $\gamma\neq 0$.
\end{corollary}
\textit{Proof.}
  For $\gamma\in\square_q$ (this includes $\gamma = 0$) the formulas were proven in Lemma
  \ref{lemma_pythagorean_triples}. So let $\gamma\in\mathbb{F}_q$ be a non-square. As the squares
  $\square_q\backslash\{ 0\}$ form a subgroup of $\mathbb{F}_q^*$,
  the non-squares have the form $\gamma\cdot \delta^2$ with $\delta\neq 0$.
  If $\alpha^2 + \beta^2 = \gamma$ then $(\alpha\delta)^2 + (\beta\delta)^2 = \gamma\delta^2$. Therefore
  the number of solutions $(\alpha,\beta)$ is the same for all non-squares and we can
  determine the number of solutions by counting: There are $q^2$ pairs $(\alpha,\beta)\in \mathbb{F}_q
  \times \mathbb{F}_q$. As there are $\frac{q-1}{2}$ squares and non-squares in $\mathbb{F}_q^*$ we obtain
  \[
    \frac{q-1}{2}\cdot\left|\mathbb{I}_\gamma\right| = q^2 - (2q-1) - \frac{q-1}{2}\cdot (q-1) =
    \frac{1}{2}q^2 - q + \frac{1}{2} = \frac{1}{2} (q-1)^2
  \]
  for $q\equiv 1\pmod 4$ and 
  \[
    \frac{q-1}{2}\cdot\left|\mathbb{I}_\gamma\right| = q^2 - 1 - \frac{q-1}{2}\cdot (q+1) = \frac{1}{2} (q-1)(q+1)
  \]
  for $q\equiv 3\pmod 4$, which yields our statement.
\hfill{$\square$}

Now we want to study the automorphism group of $\mathbb{F}_q^m$ with 
respect to $\Delta$. As already mentioned it is easy to determine the
automorphisms:

\begin{lemma}
  \label{lemma_equivalence}
  Examples of automorphisms of $\mathbb{F}_q^m$ with respect to $\Delta$ are given by:
  \begin{enumerate}
    \item $\varphi_{v}(u)=\begin{pmatrix}u_1+v_1&\dots&u_m+v_m\end{pmatrix}^T$ for $v\in\mathbb{F}_q^m$,
    \item $\widetilde{\varphi}_\alpha(u)=\alpha\cdot u$ for $\alpha\in\mathbb{F}_q^*$,
    \item $\widetilde{\varphi}_{A}(u)=A\cdot u$ for $A\in O(m,q)$, and
    \item $\widehat{\varphi}_i(u)=\begin{pmatrix}u_1^{p^i}&\dots&u_m^{p^i}\end{pmatrix}^T$
          for $i\in\{1,\dots, r-1\}$ and $q=p^r$.
  \end{enumerate}
\end{lemma}
\textit{Proof.}
  The first two cases are easy to check. For the third case we consider
  \begin{eqnarray*}
    d^2(Au,Av)&=&\langle A(u-v),A(u-v)\rangle=(u-v)^TA^TA(u-v)\\
              &=&(u-v)^T(u-v)=\langle u-v,u-v\rangle=d^2(u,v)
  \end{eqnarray*}
  and for the fourth case we have
  $$
    d^2\left(\widehat{\varphi}_i(0),\widehat{\varphi}_i(u)\right)=
    \sum_{j=1}^m \left(u_j^{p^i}\right)^2=\sum_{j=1}^m \left(u_j^2\right)^{p^i}
    =\left(\sum_{j=1}^m u_j^2\right)^{p^i}=d^2(0,u)^{p^i}.
  $$
\hfill{$\square$}

We would like to remark that the orders of the groups $O(m,q)$, $GL(m,q)$, and $\OZ(m,q)$ are as follows:
\begin{itemize}
 \item[(1)] $\left|GL(m,q)\right|=\prod\limits_{i=0}^{m-1}\left(q^m-q^i\right)$ for all $m\in\mathbb{N}$.
 \item[(2)] $\left|O(2n+1,q)\right|=2q^n\cdot\prod\limits_{i=0}^{n-1}\left(q^{2n}-q^{2i}\right)$ for $n\in \mathbb{N}$.
 \item[(3)] $\left|O(2n,q)\right|=2\left(q^n-1\right)\cdot\prod\limits_{i=1}^{n-1}\left(q^{2n}-q^{2i}\right)$
            for $n\in \mathbb{N}$ and $-1\in\square_q$.
 \item[(4)] $\left|O(2n,q)\right|=2\left(q^n+(-1)^{n+1}\right)\cdot\prod\limits_{i=1}^{n-1}\left(q^{2n}-q^{2i}\right)$
            for $n\in \mathbb{N}$ and $-1\notin\square_q$.
 \item[(5)] $|\OZ(m,q)|=\frac{q-1}{2}\cdot\left|O(m,q)\right|$
            for all $m\in\mathbb{N}\backslash\{1\}$.
\end{itemize}

As the Frobenius homomorphisms and the translations are automorphisms with respect to $\Delta$ it suffices to determine the matrix group $\Aut(m,q)\cap GL(m,q)$ of all matrices that are automorphisms with respect to $\Delta$ in order to determine the whole automorphism group. Due to Lemma \ref{lemma_equivalence} we have $\OZ(m,q)\le \Aut(m,q)\cap GL(m,q)$. Thus for dimension $m=3$ we have $(q-1)^2q(q+1)\mid\, \left|\Aut(3,q)\cap GL(3,q)\right|$. We will prove later on that $\OZ(3,q)$ is already isomorphic to $\Aut(3,q)\cap GL(3,q)$.

Firstly we summarize our knowledge on $\Aut(m,q)$:
\begin{theorem}
  \label{thm_knowledge_automorphism_group}
   We have
  \begin{itemize}
   \item[(1)] $\Aut(m,q)=A\Gamma L(m,q)$ for $2\mid q$,
   \item[(2)] $\Aut(1,q)=A\Gamma L(1,q)$,
   \item[(3)] $\Aut(2,q)\cap GL(2,q)=\OZ(2,q)$ for $2\nmid q$, $q\notin\{5,9\}$,
   \item[(4)] $\Aut(2,5)\cap GL(2,5)>\OZ(2,5)$, $\frac{\left|\Aut(2,5)\cap GL(2,5)\right|}
              {\left|\OZ(2,5)\right|}=2$, and
   \item[(5)] $\Aut(2,9)\cap GL(2,9)>\OZ(2,9)$, $\frac{\left|\Aut(2,9)\cap GL(2,9)\right|}
              {\left|\OZ(2,9)\right|}=3$.
  \end{itemize}
\end{theorem}
\textit{Proof.}
  (1) and (2) hold as for $m=1$ or $2\mid q$ all distances are integral. So in general we assume
  dimension $m\ge 2$ and odd characteristic $2\nmid q$ if not stated otherwise in the rest of this article.
  For the proof of (3), (4), and (5) we refer to \cite{integral_over_fields}.
\hfill{$\square$}

Next we prove some results on the orbits of $\mathbb{F}_q^m$ under the groups $O(m,q)$ and $\OZ(m,q)$. Therefore we need:
\begin{definition}
  By $\mathbb{P}_\tau$ we denote the set $\left\{u\in\mathbb{F}_q^m\backslash\{0\}\mid d^2(0,u)=\tau\right\}$,
  where the parameters $q$ and $m$ are provided by the context.
\end{definition}

\begin{lemma}
  \label{lemma_transitive_2}
  For every $\tau\in\mathbb{F}_q$ the group $O(2,q)$ acts transitively on $\mathbb{P}_\tau$.
\end{lemma}
\textit{Proof.}
  Firstly we consider $\tau\neq 0$. Therefore let $u=\begin{pmatrix}u_1&u_2\end{pmatrix}^T$ and
  $v=\begin{pmatrix}v_1&v_2\end{pmatrix}^T$ be two points in $\mathbb{F}_q^2$ with
  $u_1^2+u_2^2=v_1^2+v_2^2=\tau\neq 0$. With $\alpha=\frac{u_1v_1+u_2v_2}{\tau}$
  and $\beta=\frac{u_2v_1-u_1v_2}{\tau}$ we have
  $\alpha^2+\beta^2=\frac{(u_1^2+u_2^2)\cdot(v_1^2+v_2^2)}{\tau^2}=1$. Thus the matrix
  $A=\begin{pmatrix}\alpha&\beta\\-\beta&\alpha\end{pmatrix}$ is an element of $O(2,q)$ which maps
  $u$ to $v$.

  Now we deal with the remaining case $\tau=0$. If $-1\notin \square_q$ then we have $|\mathbb{P}_0|=0$. Thus we
  may assume $-1\in\square_q$. For $\alpha,\beta\in\mathbb{F}_q$ with $\alpha^2+\beta^2=0$ we have either
  $\alpha=\beta=0$ or $\alpha,\beta\neq 0$. In the latter case we have $\left(\frac{\alpha}{\beta}\right)^2=-1$,
  which has  two solutions $\frac{\alpha}{\beta}=\omega_q$ and $\frac{\alpha}{\beta}=-\omega_q$, where
  $\omega_q$ is a square root of $-1$. Thus we can write all elements of
  $\mathbb{P}_0$ as $\begin{pmatrix}\nu&\pm \nu\omega_q\end{pmatrix}^T$ with $\nu\in\mathbb{F}_q^*$. Now we apply all
  matrices of the form $B := \begin{pmatrix}\gamma&\delta\\-\delta&\gamma\end{pmatrix}$ with $\gamma^2+\delta^2=1$
  to the vector $\begin{pmatrix}1&\omega_q\end{pmatrix}^T$. By definition these matrices are elements of $O(2,q)$.
  If we parameterize $\gamma$ and $\delta$ as in Lemma \ref{lemma_pythagorean_triples} we get
  $B\begin{pmatrix}1&\omega_q\end{pmatrix}^T = \begin{pmatrix}\nu& \nu\omega_q\end{pmatrix}^T$,
  where $\nu = \frac{\tau-\tau^{-1}}{\tau+\tau^{-1}}+\frac{2}{\tau+\tau^{-1}}\cdot\omega_q$. By a small
  computation we check that
  \[
    f:\mathbb{F}_q\backslash\left\{0,\omega_q,-\omega_q\right\}\rightarrow\mathbb{F}_q\backslash\{-1,0,1\},\;\tau\mapsto
    \frac{\tau-\tau^{-1}}{\tau+\tau^{-1}}+\frac{2}{\tau+\tau^{-1}}\cdot\omega_q
  \]
  is well defined and injective. Thus all points $\begin{pmatrix}\nu& \nu\omega_q\end{pmatrix}^T$ are in the same
  orbit as $\begin{pmatrix}1&\omega_q\end{pmatrix}^T$ under the action of $O(2,q)$. As
  $\begin{pmatrix}-1&0\\0&1\end{pmatrix}\in O(2,q)$, the points $\begin{pmatrix}\nu& -\nu\omega_q\end{pmatrix}^T$
  are also contained in this orbit and the proposed statement holds.
\hfill{$\square$}

\begin{lemma}
\label{Basisergaenzungssatz}
  For each $u\in\mathbb{F}_q^3$ with $u_1^2+u_2^2+u_3^2=1$ there exist $v=\begin{pmatrix} v_1&v_2&v_3\end{pmatrix}^T$,
  $w=\begin{pmatrix}w_1&w_2&w_3\end{pmatrix}^T$ $\in\mathbb{F}_q^3$ fulfilling $v_1^2+v_2^2+v_3^2=w_1^2+w_2^2+w_3^2=1$,
  such that $\left\langle u,v\right\rangle=\left\langle u,w\right\rangle=\left\langle v,w\right\rangle=0$.
\end{lemma}
\textit{Proof.}
  The set $u^\bot$ of vectors $\tilde{v}$ solving the linear equation $\langle \tilde{v},u\rangle=0$ forms a
  2-dimensional vector space. Let $\tilde{v}\in u^\bot$ be an arbitrary element with $\tilde{v} \not= 0$.
  Then $\tilde{v}^\bot$ is also a 2-dimensional vector space and certainly $u\in \tilde{v}^\bot$. Thus we get
  $\tilde{v}^\bot \not=u^\bot$. Therefore the orthogonal vector space $u^\bot$ is non-degenerate in the sense
  of \cite[II.10.1]{huppert}. By \cite[II.10.2 b)]{huppert} there is an orthogonal basis
  $\left\{\hat{v},\hat{w}\right\}$ of $u^\bot$, i.~e.{} we have $\langle \hat{v}, \hat{w}\rangle = 0$ and
  $\langle \hat{v},\hat{v}\rangle$, $\langle \hat{w},\hat{w}\rangle$ $\neq 0$. For $\alpha, \beta\in \mathbb{F}_q$
  we obtain
  \[
    \langle \alpha \hat{v} + \beta\hat{w},\alpha \hat{v} + \beta\hat{w}\rangle = \alpha^2 \langle \hat{v},
    \hat{v}\rangle + \beta^2 \langle \hat{w}, \hat{w}\rangle.
  \]
  As $\langle \hat{v}, \hat{v}\rangle,\langle \hat{w}, \hat{w}\rangle \in \mathbb{F}_q^*$ there exist
  $\alpha, \beta\in\mathbb{F}_q$
  such that
  \[
    \langle \alpha \hat{v} + \beta \hat{w}, \alpha \hat{v} + \beta \hat{w}\rangle = 1
  \]
  by \cite[Lemma 11.1]{taylor}. Therefore there is a vector $v\in \mathbb{F}_q^3$ such that $\langle u, v\rangle
  = 0$ and $\langle v, v\rangle = 1$. Now $w$ can easily be constructed: The cross product
  $w := u \times v$ is a vector with
  $\langle u, w \rangle = \langle v, w\rangle = 0$ and
  \[
    \langle w, w \rangle = \langle u, u\rangle \cdot \langle v, v\rangle - \langle u, v\rangle^2 = 1.
  \]
\hfill{$\square$}

From the previous lemma we can easily deduce:
\begin{lemma}
  \label{lemma_transitive_5}
  The group $O(3,q)$ acts transitively on $\mathbb{P}_\tau$ for all $\tau\in\mathbb{F}_q$.
\end{lemma}
\textit{Proof.}
  For $\tau=0$ we refer to \cite[Theorem 11.6]{taylor}. Thus we may assume $\tau\neq 0$. Firstly we consider
  $\tau\in\square_q$.
  Let $\tilde{u}\in\mathbb{F}_q^3$ such that $\langle \tilde{u}, \tilde{u}\rangle = \nu^2=\tau \not= 0$.
  We put $u := \nu^{-1}\tilde{u}$. Then $\langle u,u\rangle = 1$ and by Lemma \ref{Basisergaenzungssatz} there
  are $v, w\in\mathbb{F}_q^3$ such that $A = \begin{pmatrix}u&v&w\end{pmatrix}$ is an orthogonal matrix.
  Thus the vectors $\begin{pmatrix}1&0&0\end{pmatrix}^T$ and $u$ are in the same orbit of $O(3,q)$. Thus all
  $\tilde{u}$ with $\langle \tilde{u}, \tilde{u}\rangle = \nu^2 \not= 0$ and the vectors
  $\begin{pmatrix}\pm\nu&0&0\end{pmatrix}$ are in the same orbit.

  Now we deal with the remaining cases $\tau\notin\square_q$. Let $u\in\mathbb{P}_\tau$ be an arbitrary
  vector. We show that there exists an element $A\in O(3,q)$ such that the third coordinate of $Au$
  is equal to zero. This reduces the problem to the $2$-dimensional case where we can apply Lemma
  \ref{lemma_transitive_2}, as we can extend a $2$-dimensional matrix $A'\in O(2,q)$ to a matrix
  $A\in O(3,q)$ by adding a third row and a third column consisting of a one in the diagonal and
  zeros elsewhere.

  If $u_2^2+u_3^2=\nu^2\neq 0$ then due to Lemma \ref{lemma_transitive_2} there is an element
  $A'\in O(2,q)$ which maps $\begin{pmatrix}u_2&u_3\end{pmatrix}^T$ to $\begin{pmatrix}\nu&0\end{pmatrix}^T$.
  Thus we can extend $A'$ to a desired matrix $A\in O(3,q)$ such that the third coordinate of $Au$ is equal to
  zero. Since $u_1^2+u_2^2+u_3^2\notin\square_q$ we cannot have $u_i^2+u_j^2=0$ for $i\neq j$. So we can assume
  $u_i^2+u_j^2\notin\square_q$ for $i\neq j$.

  For the remaining cases we use another technique. We set
  \[
    \mathbb{P}_{\tau,\mu}:=\Big\{v\in\mathbb{F}_q^3\backslash\{0\}\mid v_1=\mu,\,v_1^2+v_2^2+v_3^2=\tau\Big\}.
  \]
  By Lemma \ref{lemma_transitive_2} all points of $\mathbb{P}_{\tau,\mu}$ are contained in the same orbit under
  $O(3,q)$. From Corollary \ref{cor_num_sol} we deduce $\left|\mathbb{P}_{\tau,\mu}\right|=q-1$ for
  $q\equiv 1\pmod 4$ and $\left|\mathbb{P}_{\tau,\mu}\right|=q+1$ for $q\equiv 3\pmod 4$. Hence we have
  \[
    |\mathbb{P}_\tau|=\sum\limits_{\mu\in\mathbb{F}_q}\left|\mathbb{P}_{\tau,\mu}\right|
    =q\cdot\left|\mathbb{P}_{\tau,0}\right|.
  \]
  Now let us consider an arbitrary point $v\in\mathbb{P}_{\tau,\mu}$ and set
  $\lambda=v_1^2+v_2^2$. As $v_1^2+v_2^2+v_3^2=\tau\notin\square_q$ we have $\lambda\neq 0$. Due to Lemma
  \ref{lemma_transitive_2} all points $w\in\mathbb{F}_q^3$ with $w_1^2+w_2^2=\lambda$ lie in the same orbit as $v$ under
  $O(3,q)$.

  Due to Corollary \ref{cor_num_sol} we have at least $q+1$ solutions $(u_1, u_2)$ of the equation
  $u_1^2+u_2^2=\lambda$ for $q\equiv 3\pmod 4$ and $\frac{q+1}{2}$ of the $u_1$ are pairwise different. This means
  that every point in $\mathbb{P}_\tau$ lies in an orbit with at least
  $\frac{q+1}{2}\cdot\left|\mathbb{P}_{\tau,u_1}\right|=
  \frac{(q+1)^2}{2}>\frac{\left|\mathbb{P}_\tau\right|}{2}=\frac{q(q+1)}{2}$ points.
  Thus there can only be one orbit.

  For $q\equiv 1\pmod 4$ we similarly conclude that every point in $\mathbb{P}_\tau$ lies in an orbit with at least
  $\frac{(q-1)^2}{2}$ points. As $\left|\mathbb{P}_\tau\right|=(q-1)q$ and $\left|\mathbb{P}_{\tau,\mu}\right|
  =q-1$ for all $\mu\in\mathbb{F}_q$ there can exist two orbits at most and the length of every orbit has to
  be divisible by $\left|\mathbb{P}_{\tau,\mu}\right|=q-1$. If there exist exactly two orbits $\mathbb{B}_1$,
  $\mathbb{B}_2$ then we have w.l.o.g.{} $\left|\mathbb{B}_1\right|=\frac{q-1}{2}\cdot (q-1)$ and
  $\left|\mathbb{B}_2\right|=\frac{q+1}{2}\cdot(q-1)$. Due to $\left|\mathbb{B}_1\right|\,\Big|\,
  \left|O(3,q)\right|$ we have $(q-1)^2\,|\,4\cdot(q-1)q(q+1)$. Using $\gcd(q-1,q)=1$ we conclude
  $q-1\,|\,4(q+1)$. Thus we have $q-1\,|\, 8$, which is equivalent to $q\in\{3,5,9\}$. As $3\not\equiv 1\pmod 4$
  we only have to consider the cases $q= 9$ and $q = 5$. In $\mathbb{F}_9\simeq\mathbb{F}_3[x]/(x^2+1)$ we have
  $\square_9=\{0,1,2,x,2x\}$. As we have either $v_i = 0$ for some $i$ or
  $\left|\{v_1,v_2,v_3\}\cap\{1,2\}\right|\ge 2$ or $\left|\{v_1,v_2,v_3\}\cap\{x,2x\}\right|\ge 2$ there exist
  $i,j$ with $v_i^2+v_j^2\in \square_9$ in this case and we can apply our reduction to the $2$-dimensional case.

  For $q=5$, $\tau=2$ we have
  \begin{eqnarray*}
    &&\mathbb{B}_1=\left\{\begin{pmatrix}v_1&v_2&v_3\end{pmatrix}^T\mid v_1,v_2,v_3\in\{2,3\}\right\},\\
    &&\mathbb{B}_2=\left\{\begin{pmatrix}0&v_1&v_2\end{pmatrix}^T,\begin{pmatrix}v_1&0&v_2\end{pmatrix}^T,
      \begin{pmatrix}v_1&v_2&0\end{pmatrix}^T\mid v_1,v_2\in\{1,4\}\right\},
  \end{eqnarray*}
  and for $q=5$, $\tau=3$ we have
  \begin{eqnarray*}
    &&\mathbb{B}_1=\left\{\begin{pmatrix}v_1&v_2&v_3\end{pmatrix}^T\mid v_1,v_2,v_3\in\{1,4\}\right\},\\
    && \mathbb{B}_2=\left\{\begin{pmatrix}0&v_1&v_2\end{pmatrix}^T,\begin{pmatrix}v_1&0&v_2\end{pmatrix}^T,
    \begin{pmatrix}v_1&v_2&0\end{pmatrix}^T\mid v_1,v_2\in\{2,3\}\right\}.
  \end{eqnarray*}
  By considering the matrix $B=\begin{pmatrix}1&2&4\\2&1&4\\1&1&3\end{pmatrix}$ in $O(3,5)$ we conclude that
  in both cases $\mathbb{B}_1$ and $\mathbb{B}_2$ are contained in the same orbit.
\hfill{$\square$}

\begin{lemma}
  \label{lemma_one_coordinate_zero}
  For dimension $m\ge 4$ and $u\in\mathbb{F}_q^m$ there exists an element $A\in O(m,q)$
  such that the $m$-th coordinate of $Au$ is equal to zero.
\end{lemma}
\textit{Proof.}
  If one of the $u_i$ is equal to zero then there obviously exists such a matrix $A$. So we assume $u_i\neq 0$ for
  $1\le i\le m$.

  If $u_h^2+u_i^2+u_j^2=0$ for all pairwise different $1\le h,i,j\le m$ then we would have $u=0$ or $3\mid q$:
  As $m\ge 4$, there is at least one further index $k$. If we replace $u_h$ by $u_k$ then
  $u_k^2 + u_i^2 + u_j^2 = 0$ and $u_h^2 + u_i^2 + u_j^2 = 0$ results in $u_h^2 = u_k^2$. Replacing $u_i$ and
  $u_j$ by $u_k$ leads to $u_h^2 = u_i^2 = u_j^2 = u_k^2$,
  so we obtain $3u_i^2 = 0$ and thus $u = 0$ if $3\nmid q$.

  For $3\mid q$ the same computation leads to $u_i = \pm u_j$ for all $i, j$. W.l.o.g.{} let $u_1 = 1$. Then we
  have $u_i^2 + u_j^2 = 2$ for all $i, j > 1$. By Lemma \ref{lemma_transitive_2} the group $O(2,q)$
  acts transitively on $\mathbb{P}_2$. As we can extend $2$-dimensional orthogonal matrices by ones in the diagonal
  we can assume $u_i = 1$ for $1\le i \le m$. As the matrix
  \[
    A' = \begin{pmatrix} 1& 1& 1& 1\\
                         1& 1& 2& 2\\
                         1& 2& 1& 2\\
                         1& 2& 2& 1\end{pmatrix}\in O(4,q)
  \]
  maps $\begin{pmatrix} 1& 1& 1& 1\end{pmatrix}^T$ to $\begin{pmatrix} 1& 0& 0& 0\end{pmatrix}^T$, we can extend $A'$
  to a matrix $A$ in $O(m,q)$ such that $Au$ has a zero at coordinate $m$.

  So we may assume $0\neq u_{m-2}^2+u_{m-1}^2+u_m^2=:\mu$. As there exist $\alpha,\beta\in\mathbb{F}_q$ with
  $\alpha^2+\beta^2=\mu\neq 0$ by \cite[Lemma 11.1]{taylor} we can apply Lemma \ref{lemma_transitive_5} to deduce
  that there exists an element $A'\in O(3,q)$ which maps $\begin{pmatrix}u_{m-2}&u_{m-1}&u_m\end{pmatrix}^T$ onto
  $\begin{pmatrix}\alpha&\beta&0\end{pmatrix}^T$. Clearly we can extend $A'$
  to obtain the desired matrix $A\in O(m,q)$ mapping $v$ onto a point with $m$-th coordinate being equal to zero.
\hfill{$\square$}

\begin{theorem}
  \label{lemma_transitive_8}
  For dimension $m\ge 2$ the group $O(m,q)$ acts transitively on $\mathbb{P}_\tau$ for all $\tau\in\mathbb{F}_q$.
\end{theorem}
\textit{Proof.}
  We prove the theorem by induction and use Lemma \ref{lemma_transitive_2} and Lemma
  \ref{lemma_transitive_5} as induction basis. Now let $u,v\in\mathbb{P}_\tau$ be arbitrary.
  Due to Lemma \ref{lemma_one_coordinate_zero} there exist $A,B\in O(m,k)$ such that the $m$-th coordinate
  of $\tilde{u}=Au$ and the $m$-th coordinate of $\tilde{v}=Bv$ are both equal to zero. Deleting the last
  coordinate from $\tilde{u}$ and $\tilde{v}$ yields two vectors $\hat{u}$ and $\hat{v}$ in
  $\mathbb{P}_\tau$, respectively. Due to our induction hypothesis there exists an element $C'\in O(m-1,q)$
  with $C'\hat{u}=\hat{v}$. Clearly we can extend $C'$ to a matrix $C\in O(m,q)$ with $C\tilde{u}=\tilde{v}$.
  With $D=B^{-1}CA$ we have
  $D\in O(m,q)$ and $Du=v$.
\hfill{$\square$}

\begin{definition}
  \label{def_p_plus_minus}
  We set
  \[
    \mathbb{P}^+:=\bigcup_{\tau\in\square_q\backslash\{0\}}\mathbb{P}_\tau
    \quad\text{and}\quad
    \mathbb{P}^-:=\bigcup_{\tau\notin\square_q}\mathbb{P}_\tau.
  \]
\end{definition}

\begin{lemma}
  \label{lemma_orbits_1}
  For $2\nmid q$ and $m\ge 2$ the orbits of $\mathbb{F}_q^m$ under the group $\OZ(m,q)$
  are $\mathbb{P}^+$, $\mathbb{P}_0$, and $\mathbb{P}^-$.
\end{lemma}
\textit{Proof.}
  From the previous lemmas we know that $O(m,q)$ acts transitively on $\mathbb{P}_\tau$ for $2\nmid q$, $m\ge 2$,
  and $\tau\in\mathbb{F}_q$. Thus $\OZ(m,q)$ acts transitively on $\mathbb{P}^+$, $\mathbb{P}_0$, and
  $\mathbb{P}^-$. (For $A\in O(m,q)$ and $\alpha\in\mathbb{F}_q^*$ we have $B:=\alpha\cdot A\in\OZ(m,q)$ and
  $\langle Bu,Bu\rangle=\alpha^2\langle u,u\rangle$ for all $u\in\mathbb{F}_q^m$.)
\hfill{$\square$}

\begin{lemma}
  \label{lemma_orthogonal_plane}
  Let $u,v\in \mathbb{F}_q^3$ with $u, v\not= 0$. If
  $\langle u,u\rangle=\langle u,v\rangle=0$ and $\mathbb{F}_q\cdot u\neq\mathbb{F}_q\cdot v$ then we have
  $\langle v,v\rangle\in\square_q$ if $q\equiv 1\pmod 4$ and
  $\langle v,v\rangle\notin\square_q$ if $q\equiv 3\pmod 4$.
\end{lemma}
\textit{Proof.}
  If the $u_i$ are non-zero we can assume w.l.o.g.{} that $u_3=1$. From $\langle u,v\rangle=0$ we conclude
  $v_3=-v_1u_1-v_2u_2$. Using $u_1^2+u_2^2+1=0$ this results in
  \begin{eqnarray*}
    \langle v,v\rangle &=& v_1^2+v_2^2+v_1^2u_1^2+v_2^2u_2^2+2v_1v_2u_1u_2\\
    &=& -u_2^2v_1^2-u_1^2v_2^2+2v_1v_2u_1u_2\\
    &=&-(u_2v_1-u_1v_2)^2.
  \end{eqnarray*}
  As $-1\notin\square_q$ iff $q\equiv 3\pmod 4$ by Lemma \ref{lemma_root_of_minus_one} we have
  \[
  \begin{array}{ll}
  \langle v,v\rangle\notin\square_q\backslash\{0\}& \mbox{ for }q\equiv 3\pmod 4\mbox{ and}\\
  \langle v,v\rangle\in\square_q& \mbox{ for } q\equiv 1\pmod 4
  \end{array}
  \]
  in this case.

  Let us assume $q\equiv 3\pmod 4$ and $\langle v,v\rangle=0$ for a moment. As $-1\notin\square_q$ we have
  $u_1,u_2\neq 0$ using $u_1^2+u_2^2+1=0$. Thus we have $v_1=v_2\frac{u_1}{u_2}$. Inserting it yields
  $v=\begin{pmatrix}v_2\frac{u_1}{u_2}& v_2& \left(-v_2\frac{u_1}{u_2}\cdot u_1 - v_2u_2\right)\end{pmatrix}^T =
  \frac{v_2}{u_2} \cdot u$. As $u,v\neq 0$ we would have $\mathbb{F}_q\cdot v = \mathbb{F}_q\cdot u$.
  Thus we even have $\langle v,v\rangle\notin\square_q$ for $q\equiv 3\pmod 4$.

  In the remaining case we assume w.l.o.g.{} $u_3=0$. As $u_1^2+u_2^2=0$ we have $-1\in\square_q$,
  $q\equiv 1\pmod 4$, and $u_1,u_2\neq 0$. We can further assume w.l.o.g.{} $u_1=1$ and $u_2=\omega_q$, where
  $\omega_q^2=-1$. With this $\langle u,v\rangle=0$ is equivalent to $v_2=\omega_q v_1$. Thus we have
  $\langle v,v\rangle=v_3^2\in\square_q$.
\hfill{$\square$}

\begin{lemma}
  \label{lemma_orbits_2}
  For $2\nmid q$ and $m\ge 3$ the orbits of $\mathbb{F}_q^m$ under the group $\Aut(m,q) \cap GL(m,q)$ are
  $\mathbb{P}^+$, $\mathbb{P}_0$, and $\mathbb{P}^-$.
\end{lemma}
\textit{Proof.}
  As $\OZ(m,q)\le \Aut(m,q)\cap GL(m,q)$ and due to Lemma \ref{lemma_orbits_1} it may only happen that
  some elements of $\mathbb{P}^+$, $\mathbb{P}_0$, and $\mathbb{P}^-$ are contained in the same orbit. Due to
  Definition \ref{def_equivalence} $\mathbb{P}^-$ forms its own orbit. Thus only $\mathbb{P}^+$ and $\mathbb{P}_0$
  may be contained in the same orbit. Now  we show that this is not the case.

  In Section \ref{sec_graph_of_integral_distances} we introduce the graph $\mathfrak{G}_{m,q}$ of integral distances
  corresponding to $\mathbb{F}_q^m$ and its integral distances. Due to Lemma \ref{lemma_not_srg} for dimension $m=3$
  and $2\nmid q$ the graph $\mathfrak{G}_{3,q}$ is not strongly regular. Thus $\mathbb{P}^+$ and $\mathbb{P}_0$ are
  disjoint orbits.


  For $m\ge 4$ let us assume that there exists an element $u$ in $\mathbb{F}_q^m$ with $\langle u,u\rangle=0$ and
  there exists a matrix $A$ in $\Aut(m,q)\cap GL(m,q)$ with $\langle A^{-1}u,A^{-1}u\rangle\in\square_q\backslash\{0\}$.
  W.l.o.g.{} we assume $A^{-1}u=e^{(1)}$, where $e^{(i)}$ is a vector in $\mathbb{F}_q^m$ consisting of zeros and a
  single one at coordinate $i$, this is the $i$-th unit vector. So we have $Ae^{(1)}=u$ and we set
  $w^{(i)}:=A e^{(i)}\in\mathbb{F}_q^m$, $\mu_i:=\langle u,w^{(i)}\rangle$  for $2\le i\le 4$. Now we show that
  there exists a vector $v\in\mathbb{F}_q^m$ with $\langle e^{(1)},v\rangle=0$, $\langle v,v\rangle\neq 0$, and
  $\langle Ae^{(1)},Av\rangle=0$. If there exists $2\le i\le 4$ with
  $\mu_i=0$ then we may choose $v=e^{(i)}$. Otherwise we have $\mu_2,\mu_3,\mu_4\neq 0$. We remark that
  $\left(\frac{\mu_i}{\mu_j}\right)^2=-1$ is equivalent to $\left(\frac{\mu_j}{\mu_i}\right)^2=-1$ for all
  $2\le i,j\le 4$. Due to $\left(\frac{\mu_1}{\mu_2}\right)^2\cdot \left(\frac{\mu_2}{\mu_3}\right)^2\cdot
  \left(\frac{\mu_3}{\mu_1}\right)^2=1\neq -1$ there exist $i$ and $j$ with $i\neq j$,
  $\left(\frac{\mu_i}{\mu_j}\right)^2\neq -1$. We set $v:=-\mu_je^{(i)}+\mu_ie^{(j)}$ which yields
  \begin{eqnarray*}
    \langle e^{(1)},v\rangle   &=& 0,\\
    \langle v,v\rangle    &=& \mu_i^2+\mu_j^2\neq 0,\text{ and}\\
    \langle Ae^{(1)},Av\rangle &=& \langle u,-\mu_jw^{(i)}+\mu_iw^{(j)}\rangle=-\mu_j\langle u,w^{(i)}\rangle+\mu_i
    \langle u,w^{(j)}\rangle=0.
  \end{eqnarray*}
  Let $\chi$ be the characteristic function of $\square_q$, this is $\chi(\alpha)=1$ for $\alpha\in\square_q$
  and $\chi(\alpha)=0$ for $\alpha\notin\square_q$.
  We set $\tau:=\langle v,v\rangle\neq 0$ and $\nu:=\langle Av,Av\rangle$. For all
  $\lambda_1,\lambda_2\in\mathbb{F}_q$ we have
  \begin{eqnarray*}
    d^2\left(0,\lambda_1e^{(1)}+\lambda_2 v\right)    &=& \langle\lambda_1e^{(1)}+\lambda_2 v,
    \lambda_1e^{(1)}+\lambda_2 v\rangle=\lambda_1^2+\tau\cdot\lambda_2^2\text{ and}\\
    d^2\left(0,A(\lambda_1e^{(1)}+\lambda_2 v)\right) &=& \langle\lambda_1u+\lambda_2Av,\lambda_1u+\lambda_2Av\rangle=
    \nu\cdot\lambda_2^2.
  \end{eqnarray*}
  As $A\in\Aut(m,q)\cap GL(m,q)$ we have
  $\chi\left(\lambda_1^2+\tau\cdot\lambda_2^2\right)=\chi\left(\nu\cdot\lambda_2^2\right)$ for all $\lambda_1,\lambda_2$
  in $\mathbb{F}_q$. Inserting $\lambda_2=1$ yields $\chi(\nu)=\chi(\lambda_1^2+\tau)$ for all
  $\lambda_1\in\mathbb{F}_q$. Due to
  $\Big|\left\{\lambda_1^2+\tau\mid\lambda_1\in\mathbb{F}_q\right\}\Big|=\frac{q+1}{2}$ we conclude
  $\chi(\nu)=\chi\left(\lambda_1^2+\tau\right)=1$. W.l.o.g.{} we may assume $\tau=1\in\square_q\backslash\{0\}$.
  Thus for $q=p^r$ and $\alpha\in\square_q$ we have
  \[
    \nu=\chi(\alpha)=\chi(\alpha+1)=\chi((\alpha+1)+1)=\chi((\alpha+2)+1)=\dots=\chi((\alpha+p-2)+1).
  \]
  We conclude $p\,\Big|\,\left|\square_q\right|=\frac{q+1}{2}=\frac{p^r+1}{2}$, which is a contradiction.
\hfill{$\square$}

\begin{lemma}
  \label{thm_automorphismgroup}
  \[
   \Aut(3,q)\cap GL(3,q) = \OZ(3,q)\text{ for }2\nmid q.
  \]
\end{lemma}
\textit{Proof.}
  Let $A\in \Aut(3,q)\cap GL(3,q)$ be an automorphism. The idea of the proof is to use the fact
  that $A$ takes every vector $v$ of integral norm $\langle v, v\rangle \not= 0$ to another vector of integral
  norm $\not= 0$ with the aim to construct an automorphism in $\Aut(2,q)$. Using
  the classification of the $2$-dimensional automorphisms in Theorem \ref{thm_knowledge_automorphism_group}, see
  also \cite{integral_over_fields}, we conclude $A\in\OZ(3,q)$.

  Obviously, $A$ is uniquely defined by its images of $e^{(1)}=\begin{pmatrix}1&0&0\end{pmatrix}^T$,
  $e^{(2)}=\begin{pmatrix}0&1&0\end{pmatrix}^T$, and $e^{(3)}=\begin{pmatrix}0&0&1\end{pmatrix}^T$. Due to Lemma
  \ref{lemma_orbits_1} and Lemma \ref{lemma_orbits_2} we can assume $A\cdot e^{(1)}=e^{(1)}$. We set $A\cdot
  e^{(2)}=:\begin{pmatrix}\alpha&\beta&\gamma\end{pmatrix}^T$, where we have
  \[
    \alpha^2+\beta^2+\gamma^2=\nu^2\in\square_q\backslash\{0\}
  \]
  due to Lemma \ref{lemma_orbits_2}. Let $\chi:\mathbb{F}_q\rightarrow\{0,1\}$, where $\chi(\tau)=1$ iff
  $\tau\in\square_q$, be the characteristic function of $\square_q$. By applying $A$ on
  $\begin{pmatrix}\lambda&\mu&0\end{pmatrix}^T$ we obtain
  \begin{equation}
    \label{eq_chi}
    \chi\left(\lambda^2+\mu^2\right) = \chi\left(\lambda^2+2\alpha\lambda\mu+\nu^2\mu^2\right)
    \text{ for all }\lambda,\mu\in\mathbb{F}_q.
  \end{equation}
  Inserting $\lambda=-2\alpha$, $\mu=1$ yields $\chi\left(4\alpha^2+1\right) = \chi\left(\nu^2\right)=1$. Now we prove
  $\chi\left(\alpha^2+1\right) = 1$.
  Putting $\lambda = 2\tau\alpha$, $\mu = 1$ for an arbitrary $\tau\in \mathbb{F}_q$ in Equation (\ref{eq_chi})
  we obtain
  \[
    \chi\left(4\tau^2\alpha^2+1\right) = \chi\left((4\tau^2+4\tau)\alpha^2+\nu^2\right).
  \]
  Inserting $\lambda=-2(\tau+1)\alpha$, $\mu=1$ yields 
  \[
    \chi\left(4(\tau+1)^2\alpha^2+1\right)=\chi\left((4\tau^2+4\tau)\alpha^2+\nu^2\right),
  \]
  hence we get
  \[
    \chi\Big((2\tau)^2 \alpha^2 + 1\Big) = \chi\Big((2\tau+2)^2 \alpha^2 + 1\Big)
  \]
  for all $\tau\in\mathbb{F}_q$. For $\tau = 1$ we have $\chi\left(4\alpha^2+1\right) = 1$, therefore we obtain
  $\chi\left((2\tau)^2 \alpha^2 + 1\right) = 1$ for all $\tau\in \mathbb{F}_p$ (but not necessarily for all
  $\tau\in \mathbb{F}_q$). As we have $p\neq 2$, we can take $\tau = 2^{-1} \in \mathbb{F}_p$ and get
  $\chi\left(\alpha^2+1\right) = 1$.

  If we insert $\lambda=-\alpha$, $\mu=1$ in Equation (\ref{eq_chi}) we obtain
  \[
    1=\chi\left(\alpha^2+1\right)=\chi\left(\lambda^2+\mu^2\right)=\chi\left(\lambda^2+2\alpha\lambda\mu+\nu^2\mu^2
    \right)=\chi\left(\nu^2-\alpha^2\right).
  \]
  Thus a $\pi\in\mathbb{F}_q$ with $\alpha^2+\pi^2=\nu^2$ exists. Let us consider the matrix
  $B=\begin{pmatrix}1&\alpha\\0&\pi\end{pmatrix}$.
  As we have $\chi\left(\lambda^2+\mu^2\right)=\chi\Big(\lambda^2+2\alpha\lambda\mu+(\alpha^2+\pi^2)\lambda^2\Big)$
  for all $\lambda,\mu\in\mathbb{F}_q$ the matrix $B$ is an  automorphism for $\mathbb{F}_q^2$ with respect to $\Delta$.

  For $q\neq \{5,9\}$ we can apply Theorem \ref{thm_knowledge_automorphism_group}.(3) and conclude $\alpha=0$,
  $\beta^2+\gamma^2=\pi^2=\nu^2=1$. Now we set
  $A\cdot e^{(3)}=:\begin{pmatrix}\tilde{\alpha}&\tilde{\beta}&\tilde{\gamma}\end{pmatrix}^T$
  and similarly conclude $\tilde{\alpha}=0$, $\tilde{\beta}^2+\tilde{\gamma}^2=1$. Therefore $A$ has the form
  \[
    A=\begin{pmatrix}1&0&0\\0&b&\tilde{\beta}\\0&c&\tilde{\gamma}\end{pmatrix} =:
    \left(\begin{array}{rrr} 1 & 0 & 0\\0 & &\\0 & \multicolumn{2}{c}
    {\mbox{\raisebox{5pt}[-5pt]{$A'$}}}\end{array}\right).
  \]
  By applying $A$ to all vectors $\begin{pmatrix}0&\mu&\kappa\end{pmatrix}$ for $\mu,\kappa\in\mathbb{F}_q$ we see
  that $A'$ is an element of $\Aut(2,q)\cap GL(2,q)$. Due to Theorem
  \ref{thm_knowledge_automorphism_group}.(3) the matrix $A'$ is orthogonal and we conclude $A\in O(3,q)$.

  We deal with the missing cases $q\in\{5,9\}$ using the classification of the $2$-dimensional automorphism group
  $\Aut(2,q)$ as follows. Either we use the precise classification in
  \cite{integral_over_fields}  or we use an exhaustive enumeration of the elements in
  $GL(2,q)$ to conclude $\alpha=0$, $\beta^2+\gamma^2=\pi^2=\nu^2=1$ for $q=5$ and $\alpha=0$,
  $\beta^2+\gamma^2=\pi^2=\nu^2\in\{\pm 1\}$ for $q=9$.
  Now we set $A\cdot e^{(3)}=:\begin{pmatrix}\tilde{\alpha}&\tilde{\beta}&\tilde{\gamma}
  \end{pmatrix}^T$ and similarly conclude $\tilde{\alpha}=0$, $\tilde{\beta}^2+\tilde{\gamma}^2=1$ for $q=5$ and
  $\tilde{\beta}^2+\tilde{\gamma}^2 \in\{\pm 1\}$ for $q=9$. Therefore $A$ has the form
  \[
    A=\begin{pmatrix}1&0&0\\0&b&\tilde{\beta}\\0&c&\tilde{\gamma}\end{pmatrix} =:
    \left(\begin{array}{rrr} 1 & 0 & 0\\0 & &\\0 & \multicolumn{2}{c}
    {\mbox{\raisebox{5pt}[-5pt]{$A'$}}}\end{array}\right).
  \]
  Additionally we have $\langle A e^{(2)},A e^{(3)}\rangle=0$ in both cases, where we refer to
  \cite{integral_over_fields}
  or an exhaustive enumeration. Next we exclude the case $\beta^2+\gamma^2=-1$ for $q=9$. We use
  $\mathbb{F}_9\simeq\mathbb{F}_3[x]/(x^2+1)$ and assume the contrary $\beta^2+\gamma^2=-1$. As $A$ is an
  automorphism of $\mathbb{F}_9^3$ with respect to $\Delta$ we have
  \begin{eqnarray*}
    \chi\left(\lambda^2+\mu^2+\kappa^2\right)&=&\chi\Big(\lambda^2 +\left(\beta\mu+\tilde{\beta}\kappa\right)^2 +
    \left(\gamma\mu+\tilde{\gamma}\kappa\right)^2\Big)\\
    &=&\chi\Big(\lambda^2+\left(\beta^2+\gamma^2\right)\mu^2+\left(\tilde{\beta}^2+\tilde{\gamma}^2\right)\kappa^2\Big)
  \end{eqnarray*}
  for all $\lambda,\mu,\kappa\in\mathbb{F}_9$. Inserting $\lambda=1$, $\mu=1$, and $\kappa=x+2$ yields
  \[
    \chi\left(\lambda^2+\mu^2+\kappa^2\right)=\chi\left(2 + x^2 + 4x + 4\right) = \chi\left(x+2\right)=0
  \]
  and
  \[
    \chi\Big(\lambda^2+\left(\beta^2+\gamma^2\right)\mu^2+\left(\tilde{\beta}^2+\tilde{\gamma}^2\right)\kappa^2\Big)
    =\chi\Big(\left(\tilde{\beta}^2+\tilde{\gamma}^2\right)\kappa^2\Big)=1,
  \]
  a contradiction. Thus due to symmetry we have $\beta^2+\gamma^2 = \tilde{\beta}^2+\tilde{\gamma}^2 = 1$ and
  $A\in O(3,q)$ in both cases.
\hfill{$\square$}

\textit{Proof of Theorem \ref{thm_automorphismgroup_general}.}
  We prove the theorem by induction on the dimension $m$. For the induction basis we refer to Lemma
  \ref{thm_automorphismgroup}. Now let $m\ge 4$ and $A\in \Aut(m,q)\cap GL(m,q)$ be an automorphism.
  Due to Lemma \ref{lemma_orbits_1} and Lemma \ref{lemma_orbits_2} we can assume
  $A\cdot e^{(1)}=e^{(1)}$, where
  $e^{(i)}=\begin{pmatrix}0&\dots&0&\underset{i\text{-th position}}{\underbrace{1}}&0&\dots&0\end{pmatrix}^T$
  again denotes the $i$-th unit vector. For $2\le i\le m$ we set $A\cdot e^{(i)}=:
  \begin{pmatrix}v_{1,i}&\dots&v_{m,i}\end{pmatrix}^T$, where we have
  \[
    \sum\limits_{j=1}^m v_{j,i}^2=\nu_i^2\in\square_q\backslash\{0\}.
  \]
  Using a similar calculation as in the proof of Lemma \ref{thm_automorphismgroup} we obtain $v_{1,i}=0$ and
  $\sum\limits_{j=1}^m v_{j,i}^2=1$ for all $2\le i\le m$. Therefore $A$ has the form
  \[
    A=\begin{pmatrix}1&0&0&\dots\\0&v_{1,2}&v_{1,3}&\dots\\0&v_{2,2}&v_{2,3}&\dots\\\vdots&\vdots&\vdots&\ddots
    \end{pmatrix} =:
    \left(\begin{array}{rrr} 1 & 0 & \dots\\ 0 & &\\ \vdots & \multicolumn{2}{c}
    {\mbox{\raisebox{5pt}[-5pt]{$A'$}}}\end{array}\right).
  \]
  As $A$ is an automorphism of $\mathbb{F}_q^m$ with respect to $\Delta$ the matrix $A'$ is an automorphism of
  $\mathbb{F}_q^{m-1}$ with respect to $\Delta$. Due to $\sum\limits_{j=1}^m v_{j,i}^2=1$ for all $2\le i\le m$ and the
  induction hypothesis we have $A'\in O(m-1,q)$. Thus we have $A\in O(m,q)$.
\hfill{$\square$}


\section{Graph of integral distances}
\label{sec_graph_of_integral_distances}
\noindent
It turns out that it is useful to model integral point sets as cliques of certain graphs. For a given prime power $q=p^r$ and a given dimension $m$ we define a graph $\mathfrak{G}_{m,q}$ with vertex set $\mathbb{F}_q^m$, where two vertices $u$ and $v$ are adjacent if $d^2(u,v)\in\square_q$. In this section we want to study the properties of $\mathfrak{G}_{m,q}$. A motivation for this study is that the graph $\mathfrak{G}_{2,q}$ for dimension $m=2$ is a strongly regular graph. A graph is strongly regular, if there exist positive integers $k$, $\lambda$, and $\mu$ such that every vertex has $k$ neighbors, every adjacent pair of vertices has $\lambda$ common neighbors, and every nonadjacent pair has $\mu$ common neighbors, see e.~g.{} \cite{west}. If we denote the number of vertices by $v$, our graph $\mathfrak{G}_{2,q}$ has the parameters $(v,k,\lambda,\mu)=$
\[
\left(q^2,\frac{(q-1)(q+3)}{2},\frac{(q+1)(q+3)}{4}-3,\frac{(q+1)(q+3)}{4}\right)\quad\text{for }q\equiv 1\pmod 4
\]
and the parameters
\[
(v,k,\lambda,\mu)=\left(q^2,\frac{q^2-1}{2},\frac{q^2-1}{4}-1,\frac{q^2-1}{4}\right)\quad\text{for }q\equiv 3\pmod 4.
\]
See e.~g.{} \cite{michael_1} for this fact, which is easy to prove.

For $2\mid q$ or $m=1$ the graph of integral distances $\mathfrak{G}_{m,q}$ is equivalent to a complete graph on $q^m$ vertices. Thus we assume $2\nmid q$ and $m\ge 3$ in the following.

As the translations of $\mathbb{F}_q^m$ are automorphisms with respect to $\Delta$ acting transitively on the
points we know that $\mathfrak{G}_{m,q}$ is a regular graph, which means that every vertex $u$ has an equal number of neighbors, called the degree of $u$. Thus we can speak of a degree of $\mathfrak{G}_{m,q}$.

\begin{lemma}
  \label{lemma_degree}
  The degree of $\mathfrak{G}_{3,q}$ is given by\\
  \begin{align*}
    \mathcal{D}(3,q)=\begin{cases}(q-1)\cdot\frac{(q+2)(q+1)}{2} & \mbox{if }
    q\equiv 1\pmod 4, \\ (q-1)\cdot \frac{q^2+q+2}{2} & \mbox{if }q\equiv 3\pmod 4.\end{cases}
  \end{align*}
  \vspace*{-10mm}
\end{lemma}
\textit{Proof.}
  It suffices to determine the number of vectors $\begin{pmatrix}\alpha&\beta&\gamma\end{pmatrix}^T\neq0$
  fulfilling $\alpha^2+\beta^2+\gamma^2\in\square_q$. So let
  $\alpha^2+\beta^2+\gamma^2=\nu^2$. If $\nu=0$ then we have $\alpha^2+\beta^2=-\gamma^2$. Using Corollary
  \ref{cor_num_sol} we obtain $q^2-1$ solutions in this case. For $\nu\neq 0$ we have $\frac{q-1}{2}$ possible
  values for $\nu^2$. Using Corollary \ref{cor_num_sol} we
  obtain the number of solutions $(\alpha, \beta)$ of the equation $\alpha^2 + \beta^2 = \nu^2 - \gamma^2$ for
  all possible values of $\gamma$ and $\nu$. Summing up everything yields the stated formula.
\hfill{$\square$}

To determine the degree $\mathcal{D}(m,q)$ of the graph of integral distances $\mathfrak{G}_{m,q}$ in arbitrary dimension we define the three functions
\begin{eqnarray*}
  \mathcal{S}(m,q) &:=& \left|\left\{u\in\mathbb{F}_q^m\,\Big|\,\sum_{i=1}^m
  u_i^2\in\square_q\backslash\{0\}\right\}\right|,\\
  \mathcal{Z}(m,q) &:=& \left|\left\{u\in\mathbb{F}_q^m\,\Big|\,\sum_{i=1}^m u_i^2=0\right\}\right|,\text{ and}\\
  \mathcal{N}(m,q) &:=& \left|\left\{u\in\mathbb{F}_q^m\,\Big|\,\sum_{i=1}^m u_i^2\notin\square_q\right\}\right|.
\end{eqnarray*}
\noindent
We would like to remark that $\mathcal{S}(m,q)=\left|\mathbb{P}^+\right|$ and $\mathcal{N}(m,q)=\left|\mathbb{P}^-\right|$, where $\mathbb{P}^+$ and $\mathbb{P}^-$ denote the sets of Definition \ref{def_p_plus_minus}. The first few functions are given by
\begin{eqnarray*}
  \mathcal{S}(1,q) &=& q-1,\\
  \mathcal{Z}(1,q) &=& 1,  \\
  \mathcal{N}(1,q) &=& 0,\\
  \mathcal{S}(2,q) &=& \left\{\begin{array}{rcl}\frac{(q-1)^2}{2}&\text{if}&q\equiv 1\pmod
  4,\\[1ex]\frac{q^2-1}{2}&\text{if}&q\equiv 3\pmod 4,\end{array}\right.\\
  \mathcal{Z}(2,q) &=& \left\{\begin{array}{rcl}2q-1&\text{if}&q\equiv 1\pmod 4,\\1&\text{if}&q\equiv 3\pmod
  4,\end{array}\right.\text{ and}\\
  \mathcal{N}(2,q) &=& \left\{\begin{array}{rcl}\frac{(q-1)^2}{2}&\text{if}&q\equiv 1\pmod
  4,\\[1ex]\frac{q^2-1}{2}&\text{if}&q\equiv 3\pmod 4,\end{array}\right.\\
\end{eqnarray*}
\noindent
see Corollary \ref{cor_num_sol}. To determine these functions recursively we can use:
\begin{lemma}
  \label{lemma_recursion}
  Let $\mathbb{I}_0$ and $\mathbb{I}_1$ be the sets defined in Corollary \ref{cor_num_sol}. Then for
  dimension $m\ge 3$ we have
  \begin{eqnarray*}
    \mathcal{Z}(m,q) &=& \mathcal{Z}(m-2,q)\cdot\left|\mathbb{I}_0\right|+\Big(q^{m-2}-\mathcal{Z}(m-2,q)\Big)\cdot
    \left|\mathbb{I}_1\right|,\\
    \mathcal{S}(m,q) &=& \frac{q-1}{2}\cdot\Big(\mathcal{N}(m-2,q)+\mathcal{Z}(m-2,q)\Big)\cdot
    \left|\mathbb{I}_1\right|+\frac{q-3}{2}\cdot\mathcal{S}(m-2,q)\cdot\left|\mathbb{I}_1\right|\\
    &&+\mathcal{S}(m-2,q)\cdot\left|\mathbb{I}_0\right|,\\
    \mathcal{N}(m,q) &=& q^m-\mathcal{S}(m,q)-\mathcal{Z}(m,q),\text{ and}\\
    \mathcal{D}(m,q) &=& \mathcal{S}(m,q)+\mathcal{Z}(m,q)-1.
  \end{eqnarray*}
\end{lemma}
\textit{Proof.}
  We rewrite the equation $\sum\limits_{i=1}^m u_i^2=\tau$ as $u_1^2+u_2^2=\tau-\sum\limits_{i=3}^m u_i^2$ and
  apply Corollary \ref{cor_num_sol}.
\hfill{$\square$}

\begin{theorem}
  \label{thm_degree}
  Let $m\ge 1$ be arbitrary. For $q\equiv 1\pmod 4$ we have
  \begin{eqnarray*}
    \mathcal{Z}(m,q) &=& \begin{cases}q^{m-1}&\text{for }m\text{ odd},\\q^{m-1}+q^{\frac{m}{2}}-q^{\frac{m-2}{2}}&
               \text{for }m \text{ even},\end{cases}\\
    \mathcal{S}(m,q) &=& \begin{cases}\frac{1}{2}\left(q^m-q^{m-1}+q^{\frac{m+1}{2}}-q^{\frac{m-1}{2}}\right)&
    \text{for }m\text{ odd},\\\frac{1}{2}\left(q^m-q^{m-1}-q^{\frac{m}{2}}+q^{\frac{m-2}{2}}\right)&\text{for }m\text{
    even},\end{cases}\\
    \mathcal{N}(m,q) &=& \begin{cases}\frac{1}{2}\left(q^m-q^{m-1}-q^{\frac{m+1}{2}}+q^{\frac{m-1}{2}}\right)&
               \text{for }m\text{ odd},\\\frac{1}{2}\left(q^m-q^{m-1}-q^{\frac{m}{2}}+q^{\frac{m-2}{2}}\right)&
               \text{for }m\text{ even},\end{cases}\\
   \mathcal{D}(m,q) &=&
   \begin{cases}\frac{1}{2}\left(q^m+q^{m-1}+q^{\frac{m+1}{2}}-q^{\frac{m-1}{2}}\right)-1
               &\text{for }m\text{ odd},\\\frac{1}{2}\left(q^m+q^{m-1}+q^{\frac{m}{2}}-q^{\frac{m-2}{2}}\right)-1&
               \text{for }m\text{ even}.\\\end{cases}
  \end{eqnarray*}
  For $q\equiv 3\pmod 4$ we have
  \begin{eqnarray*}
    \mathcal{Z}(m,q) &=& \begin{cases}q^{m-1}&\text{for }m\text{ odd},\\ q^{m-1}+(-q)^{\frac{m}{2}} +
              (-q)^{\frac{m-2}{2}}& \text{for }m\text{ even},\end{cases}\\
    \mathcal{S}(m,q) &=&
    \begin{cases}\frac{1}{2}\left(q^m-q^{m-1}-(-q)^{\frac{m+1}{2}}-(-q)^{\frac{m-1}{2}}\right)&\text{for }m\text{ odd},\\
       \frac{1}{2}\left(q^m-q^{m-1}-(-q)^{\frac{m}{2}}-(-q)^{\frac{m-2}{2}}\right)&\text{for }m\text{ even},\end{cases}\\
    \mathcal{N}(m,q) &=& \begin{cases}\frac{1}{2}\left(q^m-q^{m-1}+(-q)^{\frac{m+1}{2}}+(-q)^{\frac{m-1}{2}}\right)&
    \text{for }m\text{ odd},\\ \frac{1}{2}\left(q^m-q^{m-1}-(-q)^{\frac{m}{2}}-(-q)^{\frac{m-2}{2}}\right)&
    \text{for }m\text{ even},\end{cases}\\
    \mathcal{D}(m,q) &=&
    \begin{cases}\frac{1}{2}\left(q^m+q^{m-1}-(-q)^{\frac{m+1}{2}}-(-q)^{\frac{m-1}{2}}\right)-1&
               \text{for }m\text{ odd},\\
               \frac{1}{2}\left(q^m+q^{m-1}+(-q)^{\frac{m}{2}}+(-q)^{\frac{m-2}{2}}\right)-1&\text{for }m
               \text{ even}.\end{cases}
  \end{eqnarray*}
\end{theorem}
\textit{Proof.}
  Induction on $m$ using Lemma \ref{lemma_recursion}.
\hfill{$\square$}

\noindent
With strongly regular graphs in mind we consider the number of common neighbors.

\begin{theorem}
  \label{thm_common_neighbors}
  If $\mathcal{A}(m,q)$ denotes the number of common neighbors of $0$ and
  $e^{(1)}=\begin{pmatrix}1&0&\dots&0\end{pmatrix}^T$ in $\mathbb{F}_q^m\backslash\{0,e^{(1)}\}$, then for
  $m\ge 1$ we have
  \[
    \mathcal{A}(m,q) = \begin{cases}\frac{q^{m-2}\cdot(q+1)^2+(-1)^{\frac{(m-1)(q-1)}{4}}\cdot q^{\frac{m-3}{2}}\cdot
                 \left(3q^2-2q-1\right)}{4}-2
                 &\text{for }m\text{ odd},\\
                 \frac{q^{m-2}\cdot(q+1)^2+2\cdot(-1)^{\frac{m(q-1)}{4}}\cdot
                 q^{\frac{m-2}{2}}\cdot\left(q-1\right)}{4}-2
                 &\text{for }m\text{ even}.\end{cases}
  \]
\end{theorem}
\textit{Proof.}
  Clearly we have $\mathcal{A}(1,q)=q-2$. For $m\ge 1$ we count the number of
  solutions $\left(v_1,\dots,v_m\right)$ of the equation system
  \begin{eqnarray*}
     v_1^2+\sum\limits_{i=2}^m v_i^2                &=& \alpha^2,\\
     \left(v_1-1\right)^2+\sum\limits_{i=2}^m v_i^2 &=& \beta^2.\\
  \end{eqnarray*}
  There are $\left(\frac{q+1}{2}\right)^2$ different pairs
  $\left(\alpha^2,\beta^2\right)$ for $\alpha,\beta\in\mathbb{F}_q$. For given
  $\alpha^2,\beta^2$ we have $v_1=\frac{\alpha^2-\beta^2+1}{2}$ and $\sum\limits_{i=2}^m
  v_i^2=\frac{4\alpha^2\beta^2-\left(\alpha^2+\beta^2-1\right)^2}{4}=-\left(\frac{\alpha^2-\beta^2-1}{2}\right)^2
  +\beta^2=:\tau$.
  Each of the $\left(\frac{q+1}{2}\right)^2$ cases leads to a specific
  $\tau\in\mathbb{F}_q$. Now let $a_\nu$ be the number of pairs $(\alpha^2, \beta^2)$ which result in $\tau=\nu$.
  Then we obtain
  \begin{equation}
    \label{eq_sum_a_i}
    \sum_{\nu\in\mathbb{F}_q} a_\nu=\left(\frac{q+1}{2}\right)^2.
  \end{equation}
  By $b_{m,\nu}$ we denote the number of vectors
  $\begin{pmatrix}v_2&\dots&v_m\end{pmatrix}\in\mathbb{F}_q^{m-1}$ with
  $\sum\limits_{j=2}^m v_j^2=\nu$. With this we have
  \begin{equation}
    \label{eq_n_c}
    \mathcal{A}(m,q) =\sum\limits_{\nu\in\mathbb{F}_q} a_\nu\cdot b_{m,\nu}-2.
  \end{equation}
  Due to $\mathcal{A}(1,q)=q-2$ and $b_{1,\nu}=0$ for $\nu\neq 0$ we have $a_0=q$.
  If $\nu,\mu\in \square_q\backslash\{0\}$ or $\nu,\mu\notin\square_q$
  then we have $b_{m,\nu}=b_{m,\mu}$. Next we show
  \begin{equation}
    \label{eq_a_s}
    a_+:=\sum_{\nu\in\square_q\backslash\{0\}} a_\nu=
    \begin{cases}\frac{(q+1)(q-1)}{8} & \text{for } q\equiv 1 \mod 4\\
    \frac{(q-1)(q-3)}{8}& \text{for } q\equiv 3 \mod 4,\end{cases}
  \end{equation}
  from which we can conclude
  \[
    a_-:=\sum_{\nu\not\in\square_q} a_\nu=
    \begin{cases}\frac{(q-1)(q-3)}{8} & \text{for } q\equiv 1 \mod 4\\
    \frac{(q+1)(q-1)}{8}& \text{for } q\equiv 3 \mod 4,\end{cases}
  \]
  due to Equation (\ref{eq_sum_a_i}). We use the information for
  dimension $m=2$. For $\nu\notin\square_q$ we have $b_{2,\nu}=0$
  and for $\nu\in\square_q\backslash\{0\}$ we have $b_{2,\nu}=2$. For $q\equiv
  3\pmod 4$ we have $b_{2,0}=1$ and for $q\equiv 1\pmod 4$
  we have $b_{2,0}=1$. Inserting this and the formula for $\mathcal{A}(2,q)$
  in Equation (\ref{eq_n_c}) yields Equation (\ref{eq_a_s}).

  Using $b_{m,0}=\mathcal{Z}(m-1,q)$, $b_{m,\nu}=\frac{2}{q-1}\cdot \mathcal{S}(m-1,q)$ for
  $\nu\in\square_q\backslash\{0\}$, and $b_{m,\nu}=
  \frac{2}{q-1}\cdot \mathcal{N}(m-1,q)$ for $\nu\notin\square_q$ we get
  \[
    \mathcal{A}(m,q)=q\cdot \mathcal{Z}(m-1,q) + a_+\cdot\frac{2}{q-1}\cdot \mathcal{S}(m-1,q) +
   a_-\cdot\frac{2}{q-1}\cdot \mathcal{N}(m-1,q)-2
  \]
  and we obtain the stated formula by using Theorem \ref{thm_degree}.
\hfill{$\square$}

So for dimension $m=3$ we have $\mathcal{A}(3,q)=\frac{q^3+5q^2-q-9}{4}$ for $q\equiv 1\pmod 4$ and $\mathcal{A}(3,q)=\frac{q^3-q^2+3q-7}{4}$ for $q\equiv 3\pmod 4$.

\begin{lemma}
  \label{lemma_not_srg}
  For odd dimension $m\ge 3$ the graph of integral distances $\mathfrak{G}_{m,q}$ is not a strongly regular graph.
\end{lemma}
\textit{Proof.}
  Let us assume that $\mathfrak{G}_{m,q}$ is strongly regular. Then there exist corresponding parameters
  $(v,k,\lambda,\mu)$ with
  \begin{eqnarray*}
    v       &=& q^m,\\
    k       &=& \mathcal{D}(m,q),\text{ and}\\
    \lambda &=& \mathcal{A}(m,q).
  \end{eqnarray*}
  For a strongly connected graph we have the identity $(v-k-1)\mu=k(k-\lambda-1)$, see e.~g.{} \cite{west}. Using
  Theorem \ref{thm_degree} and Theorem \ref{thm_common_neighbors} we can use this identity to determine $\mu$.
  For $q\equiv 1\pmod 4$ and $m\text{ odd}$ we have
  \begin{eqnarray*}
    \!\!\!\!\!\!\!\!\!\!\!\!k(k-\lambda-1)&=&\frac{q^{\frac{m-3}{2}}(q-1)(q+1)\left(q^{\frac{m-1}{2}}-1\right)
    \left(q^m+q^{m-1}+q^{\frac{m+1}{2}}-q^{\frac{m-1}{2}}-2\right)}{8},\\
    v-k-1&=&\frac{q^{\frac{m-1}{2}}\cdot(q-1)\cdot\left(q^{\frac{m-1}{2}}-1\right)}{2},\quad\text{and}\\
    \mu&=&\frac{(q+1)\cdot\left(q^m+q^{m-1}+q^{\frac{m+1}{2}}-q^{\frac{m-1}{2}}-2\right)}{4q}.
  \end{eqnarray*}
  For $q\equiv 3\pmod 4$ and $m\text{ odd}$ we obtain
  \begin{eqnarray*}
   \!\!\!\!\!\!\!\!\!\!
   k(k-\lambda-1)&=&\frac{q^{\frac{m-3}{2}}(q-1)(q+1)\left(q^{\frac{m-1}{2}}-(-1)^{\frac{m-1}{2}}\right)}{8}\cdot\\
   &&\left(q^m+q^{m-1}-(-q)^{\frac{m+1}{2}}-(-q)^{\frac{m-1}{2}}-2\right),\\
   \!\!\!\!\!\!\!\!\!\!\!\!\!\!\!\!\!\!\!\! v-k-1
 &=&\frac{q^{\frac{m-1}{2}}\cdot(q-1)\cdot\left(q^{\frac{m-1}{2}}-(-1)^{\frac{m-1}{2}}\right)}{2},\quad\text{and}\\
    \!\!\!\!\!\!\!\!\!\!\!\!\!\!\!\!\!\!\!\!\mu   &=&
 \frac{(q+1)\cdot\left(q^m+q^{m-1}-(-q)^{\frac{m+1}{2}}-(-q)^{\frac{m-1}{2}}-2\right)}{4q}.
  \end{eqnarray*}
  As for odd $m\ge 3$ the denominator of $\mu$ is divisible by $q$ and the numerator is not divisible by $q$,
  the graph of integral distances $\mathfrak{G}_{m,q}$ is not a strongly regular graph in these cases.
\hfill{$\square$}

If we accomplish the same computation for even $m$ then for $q\equiv 1\pmod 4$ we get
\begin{eqnarray*}
     \!\!\!\!\!\!\!\!\!\!\!\!
    k(k-\lambda-1)&=&\frac{q^{m-2}\cdot(q-1)\cdot(q+1)\cdot\left(q^{\frac{m}{2}}-1\right)
    \cdot\left(q^{\frac{m}{2}}+q^{\frac{m-2}{2}}+2\right)}{8},\\
    v-k-1 &=& \frac{q^{\frac{m-2}{2}}\cdot(q-1)\cdot\left(q^{\frac{m}{2}}-1\right)}{2},\quad\text{and}\\
    \mu   &=& \frac{q^{\frac{m-2}{2}}\cdot(q+1)\cdot\left(q^{\frac{m}{2}}+q^{\frac{m-2}{2}}+2\right)}
              {4},
\end{eqnarray*}
and for $q\equiv 3\pmod 4$ we obtain
\begin{eqnarray*}
    \!\!\!\!\!\!\!\!\!\!
    k(k-\lambda-1)&=&\frac{q^{m-2}(q-1)(q+1)\left(q^{\frac{m}{2}}-(-1)^{\frac{m}{2}}\right)
                     \left(q^{\frac{m}{2}}+q^{\frac{m-2}{2}}+2\cdot(-1)^{\frac{m}{2}}\right)}{8},\\
    v-k-1 &=&\frac{q^{\frac{m-2}{2}}\cdot(q-1)\cdot\left(q^{\frac{m}{2}}-(-1)^{\frac{m}{2}}\right)}{2},\quad\text{and}\\
    \mu
  &=&\frac{q^{\frac{m-2}{2}}\cdot(q+1)\cdot\left(q^{\frac{m}{2}}+q^{\frac{m-2}{2}}+2\cdot(-1)^{\frac{m}{2}}\right)}
             {4}.
\end{eqnarray*}
So in both cases we have $\mu\in\mathbb{N}$ for even dimension $m$. Therefore the graph $\mathfrak{G}_{m,q}$ could be strongly regular for even dimension $m$, and indeed this is our conjecture:
\begin{conjecture}
  \label{conj_common_neighbors}
  If $\mathcal{B}(m,q)$ denotes the number of common neighbors of $0$ and an element $v$ with $\langle v,v\rangle=0$ in
  $\mathbb{F}_q^m\backslash\{0,v\}$, then for $m\ge 2$ we have
  \[
    \mathcal{B}(m,q)=\begin{cases}\mathcal{A}(m,q)&\text{for }m\text{ even},\\
                \mathcal{A}(m,q)-(-1)^{\frac{(q-1)(m-1)}{4}}\cdot q^{\frac{m-3}{2}}\cdot\frac{q^2-1}{4}&
                \text{for }m\text{ odd}.\end{cases}
  \]
  For even dimension $m$ the graph of integral distances $\mathfrak{G}_{m,q}$ is a strongly regular graph.
\end{conjecture}
We remark that due to Lemma \ref{lemma_orbits_2} $\mathcal{B}(m,q)$ is well defined. For $q=p^1$ being a prime we have verified Conjecture \ref{conj_common_neighbors} for small values using computer calculations. More explicitly, Conjecture \ref{conj_common_neighbors} is valid for $(m=3,p\le 2029)$, $(m=4,p\le 283)$, $(m=5,p\le 97)$, $(m=6,p\le 59)$, $(m=7,p\le 31)$, and $(m=8,p\le 23)$.

\section{Maximum cardinality of integral point sets in $\mathbb{F}_q^m$}
\label{sec_maximum_cardinality}

\noindent
In Section \ref{sec_integral_point_sets} we have introduced the notion $\mathcal{I}(m,q)$ for the maximum cardinality of an integral point set over $\mathbb{F}_q^m$. As for $m=1$ or $2 \mid q$ all distances in $\mathbb{F}_q^m$ are integral, we have $\mathcal{I}(m,q)=q^m$ in these cases. We have already stated $\mathcal{I}(2,q)=q$ for $2\nmid q$ in Theorem \ref{thm_cardinality_m_2}. Combining this with the obvious bound $\mathcal{I}(m,q)\le q\cdot\mathcal{I}(m-1,q)$ we obtain
$$
  \mathcal{I}(3,q)\le q^2
$$
for $2\nmid q$ and the next open case of dimension $m=3$.

\begin{theorem}
  \label{thm_construction_hyperplane}
  If $q\equiv 1\pmod 4$ then we have $\mathcal{I}(3,q)= q^2$.
\end{theorem}
\textit{Proof.}
  Consider the point set
  \[
    \Ps:=\left\{(\alpha,\omega_q\alpha,\beta)\mid \alpha,\beta\in\mathbb{F}_q\right\}.
  \]
  This point set is an integral point set of cardinality $q^2$.
\hfill{$\square$}

We remark that the constructed point set geometrically is a hyperplane of $\mathbb{F}_q^3$.

Using the graph of integral distances $\mathfrak{G}_{m,q}$ from Section \ref{sec_graph_of_integral_distances}
the problem of determining $\mathcal{I}(m,q)$ is transferred to the well known problem
of the determination of the maximum cardinality of cliques, these are complete subgraphs, in $\mathfrak{G}_{m,q}$.
For the latter problem software packages as e.~g.{} \textsc{Cliquer} \cite{cliquer} are available.

Thus for small values $m$ and $q$ the maximum cardinality $\mathcal{I}(m,q)$ can be exactly determined
using computer calculations. In the remaining part of this section we will deal with the cases $q\nmid 2$,
$m\ge 3$. As $\mathfrak{G}_{m,q}$ consists of $q^m$ vertices one should reduce the
problem whenever possible. One possibility is to prescribe points that must be contained in the clique. Due to
Lemma \ref{lemma_orbits_2} it suffices to investigate the two cases where we prescribe $0\in \mathbb{F}_q^m$ and
an arbitrary element of $\mathbb{P}^+$ or $\mathbb{P}_0$.

Let us consider the special case of dimension $m=3$ and $q\equiv 3\pmod{4}$. Due to
$\mathcal{I}(m,q)\ge q$ we can restrict our search on cliques $\mathbb{D}$ with cardinality at
least $q+1$. Thus there exists $\gamma\in\mathbb{F}_q$ such that the hyperplane
$\left\{ \begin{pmatrix}\alpha&\beta&\gamma\end{pmatrix}^T\mid \alpha,\beta\in\mathbb{F}_q\right\}$ contains at least two points of a clique $\mathbb{D}$ with cardinality at least $q+1$. We assume $\gamma=0$ and as for $q\equiv 3\pmod{4}$ the equation $\alpha^2+\beta^2 = 0$ has the unique solution $\alpha=\beta=0$ w.l.o.g.{} we prescribe the points $\begin{pmatrix}0&0&0\end{pmatrix}^T$ and $\begin{pmatrix}1&0&0\end{pmatrix}^T$. Additionally we know the following:
Either in such a clique $\mathbb{D}$ there exists a third point in the hyperplane with third coordinate being equal to zero, or there exist two points in a hyperplane with third coordinate being equal to an element of $\mathbb{F}_q^*$, or
every hyperplane with fix third coordinate contains at least one element of the clique $\mathbb{D}$. Using theses properties we were able to determine the following values of $\mathcal{I}(3,q)$ for small $q$:

\begin{center}
   \begin{tabular}{r|rrrrrrrrrrrrr}
     $q$                & 3 &  5 & 7 & 11 &  13 &  17 & 19 & 23 & 27 & 29 & 31 &   37 &   41 \\
     \hline
     $\mathcal{I}(3,q)$ & 4 & 25 & 8 & 11 & 169 & 289 & 19 & 23 & 28 & 841 & 31 & 1369 & 1681 \\
   \end{tabular}
\end{center}


For any given point $u\in\mathbb{F}_q^3\backslash\{0\}$ the point set $u\cdot\mathbb{F}_q$ is an integral point set over $\mathbb{F}_q^3$ of cardinality $q$. For $q\equiv 3\pmod 4$ there is another nice construction of an integral point set in $\mathbb{F}_q^3$ with cardinality $q$. Firstly we construct an integral point set \textit{on a circle}, see \cite{michael_1}. Therefore we consider the field $\mathbb{F}_q':=\mathbb{F}_q[x]/(x^2+1)$. For $\zeta=\alpha+\beta x\in\mathbb{F}_q'$ with $\alpha,\beta\in\mathbb{F}_q$ we set $\overline{\zeta}:=
\alpha-\beta x\in\mathbb{F}_q'$, which mimics the complex conjugation. Now let $\zeta$ be a generator of the cyclic group $\mathbb{F}_q'\backslash\{0\}$. We define $\mathbb{D}':=\left\{\zeta\in\mathbb{F}_q'\mid \zeta\overline{\zeta}=1\right\}$. It is not difficult to check that $\mathbb{D}'$ corresponds to an integral point set over $\mathbb{F}_q^2$ of cardinality $\frac{q+1}{2}$, see \cite{michael_1}. By $\mathbb{D}$ we denote the corresponding integral point set over $\mathbb{F}_q^3$, where the third coordinates of the points are equal to zero. Now we define the set $\mathbb{L}:=\left\{\tau\mid \tau^2+1\in\square_q\right\}$ which has cardinality
$\frac{q-1}{2}$ for $q\equiv 3\pmod 4$. With this notation we can state:

\begin{lemma}
  For $q\equiv 3\mod 4$ the set $\mathbb{D}\cup\begin{pmatrix}0&0&1\end{pmatrix}\cdot \mathbb{L}$ is an
  integral point set over $\mathbb{F}_q^3$ with cardinality $q$.
\end{lemma}

\begin{lemma}
  \label{lemma_nonintegral_plane}
  For $q\equiv 3\pmod 4$ there exists a hyperplane $\mathbb{H}$ with squared distances being either $0$ or non-squares.
\end{lemma}
\textit{Proof.}
  Due to Corollary \ref{cor_num_sol} there exist $\alpha,\beta\in\mathbb{F}_q$ with $\alpha^2+\beta^2=-1$.
  We set $u:=\begin{pmatrix}\alpha&\beta&1\end{pmatrix}^T$ and
  $v:=\begin{pmatrix}-\beta&\alpha&0\end{pmatrix}^T$. This yields $\langle u,u\rangle=0$,
  $\langle v,v\rangle=-1\notin\square_q$, and $\langle u,v\rangle =0$. Now let $\mathbb{H}:=\{\tau u+\nu v\mid
  \tau,\nu\in\mathbb{F}_q\}$. The squared distance of two elements $\tau_i u+ \nu_i v\in \mathbb{H}$, $i = 1, 2$,
  is given by $\left(\nu_1 - \nu_2\right)^2\cdot \langle v,v\rangle\notin\square_q\backslash\{0\}$.
\hfill{$\square$}

\begin{corollary}
  Let $\Ps$ be an integral point set in $\mathbb{F}_q^3$ for $q\equiv 3\pmod 4$. Either $|\Ps|\le q$ or some squared
  distances are equal to zero.
\end{corollary}
\textit{Proof.}
  We consider a covering of $\mathbb{F}_q^3$ by $q$ translations of the plane of Lemma \ref{lemma_nonintegral_plane}.
\hfill{$\square$}

We remark that our two examples of integral point sets of cardinality $q$ for $q\equiv 3\pmod 4$ do not contain a squared distance being equal to zero.

As for $q\equiv 3\pmod 4$ integral point sets over $\mathbb{F}_q^3$ of cardinality $q+1$ seem to be something special we want to list the examples that we have found by our clique search. For $q=3$ we have
\[
  \Big\{(0,0,0), (1,0,0), (2,1,1), (2,2,1)\Big\}
\]
and for $q=7$ we have
\[
  \Big\{(0,0,0), (1,0,0), (0,0,1), (1,5,5), (2,1,3), (3,1,2), (5,5,1), (6,3,6)\Big\}
\]
as examples. For $\mathcal{I}(3,27)=28$ an example is given by
\begin{eqnarray*}
  \!\!\!\!\!\!\!\!\Big\{(2+2w+2w^2,2+w^2,w^2),
  (0,2w+2w^2,1+2w),
  (1,1+w+w^2,w),
  (2,0,0),\\
  \!\!\!\!\!\!\!\!(2,w^2,2+w),
  (2,2w^2,1+2w),
  (2,2w+2w^2,1+2w),
  (w,2+2w,2+2w+2w^2),\\
  \!\!\!\!\!\!\!\!(2w,2w^2,2+2w+w^2),
  (2+2w,w^2,2+w+w^2),
  (2+2w,w+2w^2,w+2w^2),\\
  \!\!\!\!\!\!\!\!(2+2w,2w+2w^2,2w),
  (w^2,2+w+w^2,1+2w^2),
  (1+w^2,2w+2w^2,2w),\\
  \!\!\!\!\!\!\!\!(0,0,0),
  (2+w^2,1+2w,2w^2),
  (1+w+w^2,w^2,2+w^2),
  (2+w+w^2,w^2,0),\\
  \!\!\!\!\!\!\!\!(1,0,0),
  (2w+w^2,1+2w+2w^2,2+2w+2w^2),
  (2+2w+w^2,2+2w^2,1),\\
  \!\!\!\!\!\!\!\!(1,0,1+w^2),
  (1+2w^2,w+w^2,2w),
  (w+2w^2,1+w^2,1+w+2w^2),\\
  \!\!\!\!\!\!\!\!(2+w+2w^2,2+w,2+w+2w^2),
  (2+w+2w^2,2+2w,2w+w^2),\\
  \!\!\!\!\!\!\!\!(1+2w+2w^2,2+w+w^2,2w+w^2),
  (1+2w+2w^2,2w+2w^2,2w)\Big\}.
\end{eqnarray*}
where we use $\mathbb{F}_{27}\simeq\mathbb{F}_3[w]/ (w^3+w^2+w+2)$.

%
%
%
%
%

For higher dimensions we know some more exact numbers, see \cite{axel_1,axel_2}:
$\mathcal{I}(4,3)=9$, $\mathcal{I}(5,3)=27$, $\mathcal{I}(6,3)=33$, $\mathcal{I}(4,5)=25$,
$\mathcal{I}(5,5)=125$, $\mathcal{I}(4,7)=49$, $\mathcal{I}(5,7)=343$, and $\mathcal{I}(4,11)=121$.

To obtain lower bounds we can consider pairs of integral point sets $\mathbb{P}_1\subset\mathbb{F}_q^{m_1}$ and
$\mathbb{P}_2\subset\mathbb{F}_q^{m_2}$, where all squared distances in $\mathbb{P}_2$ are equal to zero. An integral
point set of cardinality $\left|\mathbb{P}_1\right|\cdot\left|\mathbb{P}_2\right|$ in $\mathbb{F}_q^{m_1+m_2}$ is given by $\left\{\begin{pmatrix}u\\v\end{pmatrix}\mid u\in\mathbb{P}_1,v\in\mathbb{P}_2\right\}$.

\begin{theorem}
  \label{lemma_lower_bound_general_q_1}
  For $q\equiv 1\pmod 4$, $m\ge 1$, and $2n\le m$ we have
  \[
    \mathcal{I}(m,q)\ge q^n\cdot\mathcal{I}(m-2n,q)
    \ge q^{\left\lceil\frac{m}{2}\right\rceil},
  \]
  where we set $\mathcal{I}(0,q)=1$.
\end{theorem}
\textit{Proof.}
  We set $\mathbb{P}_2:=\left\{\begin{pmatrix}\alpha_1&\alpha_1\omega_q&\dots&\alpha_n&\alpha_n\omega_q
  \end{pmatrix}^T\mid \alpha_1,\dots,\alpha_n\in\mathbb{F}_q\right\}$ in the above described construction.
  Thus we have $\mathcal{I}(m,q)\ge q^n\cdot\mathcal{I}(m-2n,q)$ for all $2n\le m$. The remaining inequality
  can be proven by induction on $m$ using $\mathcal{I}(m,q)\ge q$.
\hfill{$\square$}

We would like to remark that the lower bound of Theorem \ref{lemma_lower_bound_general_q_1} is sharp for $m\le 3$
and $q=5$, $m\le 5$.


\begin{lemma}
  \label{lemma_all_zero_m_4}
  There exists an integral point set $\mathbb{P}_2$ in $\mathbb{F}_q^4$ of cardinality $q^2$, where all
  squared distances are equal to zero.
\end{lemma}
\textit{Proof.}
  Let $(\alpha,\beta)$ be a solution of $\alpha^2+\beta^2=-2$ in $\mathbb{F}_q$. By Corollary \ref{cor_num_sol}
  there are at least $q-1\ge 1$ such solutions. We consider the vectors
  $u=\begin{pmatrix}\alpha&\beta&1&1\end{pmatrix}^T$ and $v=\begin{pmatrix}-\beta&\alpha&-1&1\end{pmatrix}^T$.
  Obviously $u$ and $v$ are linearly independent and fulfill $\langle u,v\rangle=0$, $\langle u,u\rangle=0$,
  and $\langle v,v\rangle=0$. We set $\mathbb{P}_2:=\left\{
  \tau u+\nu v \mid \tau,\nu\in\mathbb{F}_q\right\}$. It suffices to check $d^2\left(0,\tau u+\nu v\right)=0$ for all
  $\tau,\nu\in\mathbb{F}_q$. Indeed we have
  \[
    d^2\left(0,\tau u+\nu v\right)=\langle \tau u+\nu v,\tau u+\nu v\rangle=\tau^2\langle u,u\rangle +2\tau\nu
   \langle u,v\rangle+\nu^2\langle v,v\rangle=0.
  \]
\hfill{$\square$}

\begin{theorem}
  \label{lemma_lower_bound_general_q_2}
  For $q\equiv 3\pmod 4$, $m\ge 1$, and $4n\le m$ we have
  \[
    \mathcal{I}(m,q)\ge q^{2n}\cdot\mathcal{I}(m-4n,q)\ge
     q^{2\cdot\left\lfloor\frac{m}{4}\right\rfloor+\Big\lceil\frac{m}{4}-\left\lfloor\frac{m}{4}\right\rfloor\Big\rceil}
    \ge q^{\left\lfloor\frac{m}{2}\right\rfloor},
  \]
  where we set $\mathcal{I}(0,q)=1$.
\end{theorem}
\textit{Proof.}
  We choose $\mathbb{P}_2$ as the $n$-fold cartesian product from the integral point set of Lemma
  \ref{lemma_all_zero_m_4} in the construction described above Lemma \ref{lemma_lower_bound_general_q_1}.
  Thus we have $\mathcal{I}(m,q)\ge q^{2n}\cdot\mathcal{I}(m-4n,q)$ for all $2n\le m$.
  The remaining inequality can by proven by induction on $m$ using $\mathcal{I}(m,q)\ge q$.
\hfill{$\square$}

The lower bound of Theorem \ref{lemma_lower_bound_general_q_2} is sharp for $m\le 2$
and $q=7,11$, $m=3,4$.


\section{Conclusion and outlook}
\label{sec_conclusion}

\noindent
For the study of discrete structures the knowledge of their automorphism group is very important. In Section \ref{sec_automorphism_group} we have completed the determination of the automorphism group of $\mathbb{F}_q^m$ with respect to integral distances.

The graphs $\mathfrak{G}_{q,m}$ of integral distances are interesting combinatorial objects. We were able to determine a few parameters and properties, but the large part remains unsettled. It would be nice to have a proof of Conjecture \ref{conj_common_neighbors}, which maybe is not too difficult. We would like to remark that for all dimensions $m\ge 3$ the graph of integral distances $\mathfrak{G}_{m,q}$ is at least a slight generalization of a strongly regular graph, a so-called three class association scheme.

Section \ref{sec_maximum_cardinality} provides a first glimpse on the maximum cardinalities $\mathcal{I}(m,q)$ of integral point sets over $\mathbb{F}_q^m$. It remains a task for the future to determine some more exact numbers or lower and upper bounds. For small $q$ we have no idea for a general construction of integral point sets with maximum cardinality. A detailed analysis of the parameters of the $3$-class association schemes including the eigenvalues of the corresponding graphs could be very useful to use some general upper bounds on clique sizes. A geometrical description of the point sets achieving $\mathcal{I}(3,q)=q+1$ for $q\equiv 3\pmod 4$ would be interesting.

There are some similarities between integral point sets over $\mathbb{F}_q^m$ and integral point sets over Euclidean spaces $\mathbb{E}^m$. For example the constructions which lead to the maximum cardinality $\mathcal{I}(m,q)$ in $\mathbb{F}_q^m$ often coincide with the constructions which lead to integral point sets over $\mathbb{E}^m$ with minimum diameter, see \cite{phd_kurz,paper_laue,sascha-alfred}.



\begin{thebibliography}{10}

\bibitem{algo}
A.~Antonov and M.~Brancheva, \emph{Algorithm for finding maximal {D}iophantine
  figures}, Spring Conference 2007 of the Union of Bulgarian Mathematicians,
  2007.

\bibitem{Blokhuis-1984}
A.~Blokhuis, \emph{On subsets of ${G}{F}(q\sp 2)$ with square differences},
  Indag. Math. \textbf{46} (1984), 369--372.

\bibitem{0561.12009}
\bysame, \emph{On subsets of $gf(q\sp 2)$ with square differences}, Indag.
  Math. \textbf{46} (1984), 369--372.

\bibitem{0945.51002}
A.~Blokhuis, S.~Ball, A.~E. Brouwer, L.~Storme, and T.~Sz\H{o}nyi, \emph{On the
  number of slopes of the graph of a function defined on a finite field}, J.
  Comb. Theory, Ser. A \textbf{86} (1999), no.~1, 187--196.

\bibitem{1086.52001}
P.~Brass, W.~Moser, and J.~Pach, \emph{Research problems in discrete geometry},
  Springer, 2005.

\bibitem{Dimiev-Setting}
S.~Dimiev, \emph{A setting for a {D}iophantine distance geometry}, Tensor
  (N.S.) \textbf{66} (2005), no.~3, 275--283. \MR{MR2189847}

\bibitem{banach}
R.~E. Fullerton, \emph{Integral distances in banach spaces}, Bull. Amer. Math.
  Soc. \textbf{55} (1949), 901--905.

\bibitem{UPIN}
R.K. Guy, \emph{Unsolved problems in number theory. 2nd ed.}, Unsolved Problems
  in Intuitive Mathematics. 1. New York, NY: Springer- Verlag. 285 p., 1994.

\bibitem{integral_distances_in_point_sets}
H.~Harborth, \emph{Integral distances in point sets}, Butzer, P. L. (ed.) et
  al., Karl der Grosse und sein Nachwirken. 1200 Jahre Kultur und Wissenschaft
  in Europa. Band 2: Mathematisches Wissen. Turnhout: Brepols, 1998,
  pp.~213--224.

\bibitem{huppert}
B.~Huppert, \emph{Endliche {G}ruppen. i}, Die Grundlehren der mathematischen
  Wissenschaften in Einzeldarstellungen. 134. Berlin-Heidelberg-New York:
  Springer-Verlag. 793 p., 1967.

\bibitem{michael_1}
M.~Kiermaier and S.~Kurz, \emph{Maximal integral point sets in affine planes
  over finite fields}, Discrete Math. \textbf{309} (2009), no.~13, 4564--4575.

\bibitem{1062.00002}
M.~Kleber, \emph{Encounter at far point}, Math. Intell. \textbf{30} (2008),
  no.~1, 50--53.

\bibitem{axel_1}
A.~Kohnert and S.~Kurz, \emph{Integral point sets over $\mathbb{Z}_n^m$},
  Electronic Notes in Discrete Mathematics \textbf{27} (2006), 65--66.

\bibitem{axel_2}
\bysame, \emph{Integral point sets over $\mathbb{Z}_n^m$}, Discrete Applied
  Mathematics \textbf{157} (2009), no.~9, 2105--2117.

\bibitem{kreisel}
T.~Kreisel and S.~Kurz, \emph{There are integral heptagons, no three points on
  a line, no four on a circle}, Discrete Comput. Geom. \textbf{39} (2008),
  786--790.

\bibitem{phd_kurz}
S.~Kurz, \emph{Konstruktion und {E}igenschaften ganzzahliger {P}unktmengen},
  Ph.D. thesis, Bayreuth. Math. Schr. 76. Universit\"at Bayreuth, 2006.

\bibitem{integral_over_fields}
\bysame, \emph{Integral point sets over finite fields}, Australas. J. Comb.
  \textbf{43} (2009), 3--29.

\bibitem{paper_laue}
S.~Kurz and R.~Laue, \emph{Upper bounds for integral point sets}, Australas. J.
  Comb. \textbf{39} (2007), 233--240.

\bibitem{sascha-alfred}
S.~Kurz and A.~Wassermann, \emph{On the minimum diameter of plane integral
  point sets}, Ars Combin. \textbf{101} (2011), 265--287.

\bibitem{cliquer}
S.~Niskanen and P.~R.~J. \"Osterg{\aa}rd, \emph{Cliquer user's guide, version
  1.0}, Tech. Report T48, Communications Laboratory, Helsinki University of
  Technology, Espoo, Finland, 2003.

\bibitem{0229.12019}
L.~R\'edei, \emph{L\"uckenhafte polynome \"uber endlichen {K}\"orpern. (gap
  polynomials over finite fields)}, Mathematische Reihe. Bd. 42.
  Basel-Stuttgart: Birkh\"auser Verlag 270 S., 1970.

\bibitem{1117.05018}
P.~Sziklai, \emph{Directions in $\text{AG}(3,p)$ and their applications}, Note
  Mat. \textbf{26} (2006), no.~1, 121--130.

\bibitem{taylor}
D.~E. Taylor, \emph{The geometry of the classical groups}, Sigma Series in Pure
  Mathematics Volume 9. Heldermann Verlag Berlin, 1992.

\bibitem{west}
D.~B. West, \emph{Introduction to graph theory}, 2nd ed., New Delhi:
  Prentice-Hall of India. 608~p., 2005.

\end{thebibliography}

\end{document}